\documentclass[12pt]{amsart}

\usepackage{color,graphicx}
\usepackage{psfrag}
\usepackage{amssymb}                    
\usepackage{amsmath}                    
\usepackage{amsthm}                     
\usepackage{latexsym}                   
\usepackage{subfigure}
\usepackage{epsfig}                     
\usepackage{pstricks,pst-node,pst-text,pst-tree}
\usepackage{amscd}
\usepackage{amsmath}
\usepackage{psfrag}
\usepackage{xypic}

\newtheorem{theorem}{Theorem}[section]
\newtheorem{corollary}[theorem]{Corollary}
\newtheorem{lemma}[theorem]{Lemma}
\newtheorem{conjecture}[theorem]{Conjecture}
\newtheorem{proposition}[theorem]{Proposition}

\theoremstyle{definition}
\newtheorem{definition}[theorem]{Definition}

 \newtheorem{example}[theorem]{Example}
 \newtheorem{remark}[theorem]{Remark}
\newtheorem{algorithm}[theorem]{Algorithm}

\DeclareMathOperator{\Spec}{Spec}
\DeclareMathOperator{\Proj}{Proj}

\DeclareMathOperator{\inn}{in}
\DeclareMathOperator{\pos}{pos}

\DeclareMathOperator{\negg}{neg}
\DeclareMathOperator{\supp}{supp}

 \newcommand{\GL}{\ensuremath{\operatorname{GL}}}
 \newcommand{\Hom}{\ensuremath{\operatorname{Hom}}}
 
 \newcommand{\K}{\Bbbk}

 \newcommand{\N}{\ensuremath{\mathbb{N}}}

 \newcommand{\Q}{\ensuremath{\mathbb{Q}}}
  
 \newcommand{\SL}{\ensuremath{\operatorname{SL}}}

 \newcommand{\Z}{\ensuremath{\mathbb{Z}}}

 \newcommand{\conv}{\ensuremath{\operatorname{conv}}}
 \newcommand{\diag}{\ensuremath{\operatorname{diag}}}
 
 \newcommand{\ghilb}{\ensuremath{G}\operatorname{-Hilb}}
 \newcommand{\git}{\ensuremath{\operatorname{/\!\!/}}}
 
  \newcommand{\hilbg}{\operatorname{Hilb}^{\ensuremath{G}}}

 \newcommand{\relint}{\ensuremath{\operatorname{relint}}}
 \newcommand{\st}{\ensuremath{\operatorname{\bigm{|}}}}

 \numberwithin{equation}{section}
\begin{document}

 \bibliographystyle{plain} 
 
 \title[McKay quiver representations II: Gr\"{o}bner
 techniques]{Moduli of McKay quiver representations II: Gr\"{o}bner
   basis techniques}
 
 \author{Alastair Craw} \address{Department of Mathematics, University of Glasgow, Glasgow G12 8QW}
 \email{craw@maths.gla.ac.uk}
 
 \author{Diane Maclagan} \address{Department of Mathematics,
       Hill Center-Busch Campus,
       Rutgers University,
       110 Frelinghuysen Rd,
       Piscataway, NJ 08854} \email{maclagan@math.rutgers.edu}

 \author{Rekha R.\ Thomas} \address{Department of Mathematics, University
   of Washington, Seattle, WA
   98195}\email{thomas@math.washington.edu}

\date{\today}

 \begin{abstract}
   In this paper we introduce several computational techniques for the
   study of moduli spaces of McKay quiver representations, making use
   of Gr\"obner bases and toric geometry.  For a finite abelian group
   \(G\subset \GL(n,\K)\), let $Y_{\theta}$ be the coherent component
   of the moduli space of $\theta$-stable representations of the McKay
   quiver.  Our two main results are as follows: we provide a simple
   description of the quiver representations corresponding to the
   torus orbits of $Y_\theta$, and, in the case where $Y_\theta$
   equals Nakamura's $G$-Hilbert scheme, we present explicit equations for a cover by local coordinate
charts.  The latter theorem corrects the first
   result from Nakamura~\cite{Nakamura}.  The techniques introduced
   here allow experimentation in this subject and give concrete
   algorithmic tools to tackle further open questions.  To illustrate
   this point, we present an example of a nonnormal $G$-Hilbert
   scheme, thereby answering a question raised by Nakamura.
 \end{abstract}

 \maketitle

 \section{Introduction}
 For a finite subgroup $G\subset\SL(2,\mathbb C)$, McKay~\cite{Mckay}
 observed a connection between the representation theory of $G$, as
 encoded in the {\em McKay quiver}, and the geometry of the minimal
 resolution of $\mathbb C^2/G$.  This connection was made explicit in
 the work of Kronheimer~\cite{Kronheimer} and
 Ito--Nakamura~\cite{ItoNakamura}, who showed that the moduli spaces
 $\mathcal{M}_{\theta}$ of $\theta$-stable representations of the
 McKay quiver are isomorphic to the minimal resolution of $\mathbb
 C^2/G$.  Bridgeland-King-Reid~\cite{BKR} subsequently proved that for
 finite subgroups $G\subset \SL(3,\mathbb C)$, each of the moduli
 spaces $\mathcal{M}_{\theta}$ is isomorphic to some projective
 crepant resolution of $\mathbb C^3/G$. While the paper \cite{BKR}
 describes only the special case where $\mathcal{M}_{\theta}$ is the
 $G$-Hilbert scheme $\ghilb$, the method extends to the moduli spaces
 $\mathcal{M}_{\theta}$ for any generic parameter $\theta$ (see
 Craw-Ishii~\cite{CrawIshii}). 
 For a finite subgroup $G \subseteq \SL(n,\mathbb C)$ with $n \geq 4$,
 or for a finite subgroup $G \subset \GL(n,\mathbb C)$ with $n \geq
 3$, the moduli spaces $\mathcal{M}_{\theta}$ are no longer
 necessarily irreducible.  In Craw-Maclagan-Thomas~\cite{CMT1}, we
 introduced for finite abelian subgroups $G\subset \GL(n,\K)$, an explicit
 construction of an irreducible component $Y_{\theta}$ of
 $\mathcal{M}_{\theta}$ that is birational to $\mathbb A_{\K}^n/G$; we
 call this component the \emph{coherent component} of the moduli space
 $\mathcal{M}_{\theta}$.  

 This paper introduces several computational techniques for the study
 of moduli spaces of McKay quiver representations, making use of
 Gr\"obner bases and toric geometry.  By studying $\theta$-stable
 quiver representations in their equivalent guise as
 \emph{$G$-constellations}, we are able to study the corresponding
 modules using Gr\"{o}bner theory. Our first main result determines
 whether a given $\theta$-stable $G$-constellation corresponds to a
 point on the coherent component $Y_\theta$. In addition, when
 $\mathcal{M}_\theta \cong \ghilb$ we provide an explicit description
 of local coordinate charts on the original, irreducible version of
 the $G$-Hilbert scheme, $\hilbg$, introduced by
 Nakamura~\cite{Nakamura}.  The techniques introduced here give
 concrete algorithmic tools to tackle further open questions.  To
 illustrate this point we answer the question raised by Nakamura as to
 whether $\hilbg$ is normal by exhibiting a subgroup
 \(G\subset\GL(6,\K)\) for which $\hilbg$ is not normal.  Thus
 $\hilbg$ is an example of a nonnormal toric variety arising naturally
 in a geometric context.

 Note that Sardo Infirri~\cite{SI2} studied the moduli spaces
 \(\mathcal{M}_\theta\) for a finite abelian subgroup \(G\subset
 \GL(n,\K)\), and claimed that each \(\mathcal{M}_\theta\) was a toric variety.
 Examples~\ref{ex:reducible} and \ref{ex:nonnormal} provide
 counterexamples to this statement.

 We now describe the results in more detail.  Let $G\subset \GL(n,\K)$
 be a finite abelian subgroup, let $S = \K[x_1,\dots,x_n]$, and write
 $A:= \oplus_{\rho} S {\bf e}_\rho$ for the \(G\)-equivariant
 \(S\)-module with one generator for each irreducible representation
 \(\rho\) of $G$.  The \emph{McKay module} of \(G\) is the
 \(A\)-module
 \[
 M_G:= \langle x_i \mathbf{e}_{\rho} -
 \mathbf{e}_{\rho \rho_i}\in A : 1 \leq i \leq n, \rho \text{
 irreducible} \rangle.
 \]
 Each vector ${\bf w}\in (\mathbb Q^n_{\geq 0})^*$ determines the
 slice \(P^\vee_{\bf w} := \{\mathbf{v} \in (\Q^{r})^* : w_i+v_\rho -
 v_{\rho\rho_i}\geq 0\}\) of a polyhedral cone \(P^\vee\) that arises
 naturally from the geometric invariant theory construction of
 $Y_\theta$ (see Section~\ref{sec:moduli}).  For a parameter $\theta$
 in the GIT parameter space $\Theta$ (see Definition~\ref{defn:Theta}),
 the vector ${\bf w}\in (\mathbb Q^n_{\geq 0})^*$ determines a unique
 distinguished point of \(Y_\theta\), and hence a \emph{distinguished
   \(\theta\)-semistable \(G\)-constellation} which we denote
 \(A/M_{\theta,{\bf w}}\). The following result is proved in
 Theorem~\ref{thm:main}.

 \begin{theorem}
 \label{thm:mainintro}
 For $\theta \in \Theta$ and $\mathbf{w} \in (\mathbb Q^n_{\geq
   0})^*$, let $\mathbf{v} \in P^{\vee}_{\mathbf{w}}$ be any vector
 satisfying $\theta \cdot \mathbf{v} \leq \theta \cdot \mathbf{v}'$
 for all $\mathbf{v}' \in P^{\vee}_{\mathbf{w}}$. Then the following
 $G$-constellations coincide:
 \begin{enumerate}
 \item The distinguished \(\theta\)-semistable $G$-constellation
   \(A/M_{\theta,{\bf w}}\);
 \item The cyclic \(A\)-module $A/M_b$, where $M_b\subset A$ is the
   left $A$-ideal generated by $\{x_i\mathbf{e}_{\rho
   }-b_i^{\rho}\mathbf{e}_{\rho \rho_i} : 1 \leq i \leq n, \rho \in
   G^* \}$, and $b = (b_i^\rho)$ satisfies
 \[
 b_i^\rho =
 \left\{\begin{array}{cl} 1 & \text{if } w_i+v_\rho - v_{\rho\rho_i} =
 0 \\ 0 & \text{if } w_i+v_\rho - v_{\rho\rho_i} > 0
 \end{array}\right.;
 \] 
\item The cyclic \(A\)-module $A/\inn_{(\mathbf{v},\mathbf{w})}(M_G)$,
  where $\inn_{({\bf v},{\bf w})}(M_G)$ is the {\em initial
    module} of $M_G$ with respect to $({\bf v},{\bf w})$.
 \end{enumerate}
 \end{theorem}

 Theorem~\ref{thm:mainintro} provides a simple algorithm for computing
 $G$-constellations.  The algorithm requires that one solves a linear
 program and then calculates an initial module.  The first task is
 straightforward, and the latter is particularly simple here.
 
 Theorem~\ref{thm:mainintro} can be simplified in the case where
 \(\mathcal{M}_\theta \cong \ghilb\) and $Y_\theta \cong \hilbg$ as
 follows.  The inclusion of $G$ into $(\K^*)^n$ gives a map
 \(\deg\colon \Z^n \rightarrow \Hom(G, \K^*)\) whose kernel \(M\) is a
 lattice. Write $I_M:= \langle x^{\mathbf{u}} -x^{\mathbf{u}'} :
 \mathbf{u}, \mathbf{u}' \in \mathbb N^n, \mathbf{u}-\mathbf{u}' \in M
 \rangle$ for the lattice ideal and $\inn_{\bf w}(I_M)$ for the {\em
   initial ideal} of $M$ with respect to ${\bf w} \in (\mathbb
 Q^n_{\geq 0})^*$.  The following result is proved in
 Proposition~\ref{prop:vertexideal}.

 \begin{corollary}
  \label{coro:mainintro}
  Let $J \subseteq S$ be a monomial ideal defining a \(G\)-cluster
  \([J]\in \ghilb\).  Then \([J]\) lies in the coherent component
  \(\hilbg\) if and only if $J=\inn_{\mathbf{w}}(I_M)$ for some
  $\mathbf{w} \in (\mathbb Q^n_{\geq 0})^*$.
 \end{corollary} 

 Using Corollary~\ref{coro:mainintro}, we exhibit a finite subgroup
 \(G\subset \GL(3,\K)\) and a monomial ideal \(J\subseteq S\) such
 that \([J]\in \ghilb\) does not lie on the coherent component
 \(\hilbg\). This provides a counterexample to the statements of
 Nakamura~\cite[Corollary~2.4, Theorem~2.11]{Nakamura}, where it is
 claimed that \emph{every} monomial ideal defining a \(G\)-cluster
 gives a point of the component \(\hilbg\) (the main result of that
 paper, that $\hilbg$ is a crepant resolution of $\mathbb A^3_{\K}/G$
 when $G \subset \SL(3, \K)$, is nevertheless correct).

Stillman--Sturmfels--Thomas~\cite{AlgoTHS} established that all
 monomial ideals in the coherent component of the toric Hilbert scheme
 (see \cite{PeevaStillman}) are initial ideals of an associated toric ideal.
 Haiman--Sturmfels~\cite{HaimanSturmfels} generalized the definition
 of toric Hilbert schemes in their work on multigraded Hilbert
 schemes, including \(\ghilb\) as a special case, so
 Corollary~\ref{coro:mainintro} extends the result of \cite{AlgoTHS} to this
 case.  See Ito~\cite{Ito} for details in the \(\ghilb\) context for
 finite abelian \(G\subset \GL(2,\K)\).

 Our second main result constructs a cover of \(\hilbg\) by local
 coordinate charts.  Just as
 Corollary~\ref{coro:mainintro} is the appropriate refinement of
 \cite[Corollary~2.4]{Nakamura}, the following result (presented in
 Theorem~\ref{thm:coverhilbg}) provides the correct statement in place
 of Nakamura~\cite[Theorem~2.11]{Nakamura}.

 \begin{theorem}
 \label{thm:secondintro}
 The scheme \(\hilbg\) is covered by affine charts \(\Spec \K[A_J]\)
 indexed by monomial ideals \(J=\inn_{\mathbf{w}}(I_M)\) for
 $\mathbf{w} \in (\mathbb Q^n_{\geq 0})^*$, where \(A_J\) is a
 semigroup associated to \(J\).
 \end{theorem}

 Theorem~\ref{thm:secondintro} enables us to present the universal
 $G$-cluster over $\Spec \K[A_J]$ in an economical way (see
 Corollary 5.5). In addition, we exhibit a finite subgroup of
 \(\GL(6,\K)\) and an ideal \(J = \inn_{\mathbf{w}}(I_M)\) for which
 \(\Spec \K[A_J]\) is not normal (see Example~\ref{ex:nonnormal} and
 Corollary~\ref{coro:hilbgnormal}).  This answers the question raised
 by Nakamura~\cite[Remark~2.10]{Nakamura} as to whether $\hilbg$ is
 normal.

 \begin{corollary}
 \label{coro:secondintro}
 Nakamura's \(G\)-Hilbert scheme \(\hilbg\) is not normal in general.
 \end{corollary}
 
 While there is an extensive literature on nonnormal toric varieties
 (see \cite{GBCP}), the focus has been on applications such as integer
 programming (see, for example, \cite{HostenThomas}).  On the other
 hand, the standard definition of a toric variety in algebraic
 geometry assumes normality.  Corollary~\ref{coro:secondintro}
 therefore provides an example of a nonnormal toric variety 
 arising naturally in algebraic geometry.

 We now explain the division into sections.  Section~\ref{sec:moduli}
 reviews the construction of the moduli spaces $\mathcal{M}_{\theta}$,
 and recalls the main result from \cite{CMT1}.
 Section~\ref{sec:grobner} reviews some well-known facts from the
 theory of Gr\"{o}bner bases, and gives our first Gr\"{o}bner bases
 result for \(G\)-constellations. In
 Section~\ref{sec:distinguishedGcons} we establish
 Theorem~\ref{thm:mainintro}, and Corollary~\ref{coro:mainintro}.
 Finally, in Section~\ref{sec:nonnormal} we prove
 Theorem~\ref{thm:secondintro} and Corollary~\ref{coro:secondintro}.

 \medskip

 \noindent \textbf{Conventions.} For an integer matrix \(C\), let $\N
 C$ denote the semigroup generated by the columns of \(C\). Similarly,
 $\Z C$ denotes the lattice, $\Q_{\geq 0} C$ the rational cone and
 \(\Q C\) the rational vector space generated by columns of \(C\).
 For $\mathbf{u}=(u_1,\dots,u_m), \mathbf{u}'=(u'_1,\dots,u'_m) \in
 \mathbb N^m$ we write $\mathbf{u} \leq \mathbf{u}'$ if $u_i \leq
 u'_i$ for $1 \leq i \leq m$.  By a point of a scheme over \(\K\) we
 mean a closed point.  We write \(\K^*\) for the one-dimensional
 algebraic torus.  

 \medskip

 \noindent \textbf{Acknowledgements.} We would like to thank Bernd
 Sturmfels for bringing us together.  The original observation of a
 link between $\ghilb$ and the toric Hilbert scheme is due to him.  We
 also thank Iain Gordon, Mark Haiman, Akira Ishii, S.~Paul Smith and
 Bal\'{a}zs Szendr\H{o}i for useful comments and discussions.
 Finally, we thank the organizers of PCMI 2004 for providing a
 stimulating environment where part of this paper was written.  The
 second and third authors were partially supported by NSF grants
 DMS-0500386 and DMS-04010147 respectively.

 \section{McKay quiver representations and
  \protect$G\protect$-constellations}
 \label{sec:moduli} 
 We review the construction of the moduli spaces of McKay quiver
 representations for a finite abelian subgroup \(G\subset\GL(n,\K)\)
 of order \(r\), where \(\K\) is an algebraically closed field whose
 characteristic does not divide \(r\).  See \cite{CMT1} for a more
 leisurely introduction.  We also recall the equivalent
 module-theoretic formulation of McKay quiver representations, where
 they are known as \(G\)-constellations.

 \subsection{Moduli of McKay quiver representations}
 Since $G$ is abelian, we may assume that $G$ is contained in the
 subgroup $(\K^*)^n$ of diagonal matrices with nonzero entries in
 $\GL(n,\K)$.  We thus get $n$ elements $\rho_1,\dots, \rho_n$ of the
 dual group of characters $G^*:=\Hom(G,\K^*)$, defined by setting
 $\rho_i(g)$ to be the $i$th diagonal element of the matrix for $g$.
 The elements $\rho_1,\dots, \rho_n$ generate the group $G^*$.

 \begin{definition}
 \label{defn:quiverrep}
 The \emph{McKay quiver} of $G\subset \GL(n,\K)$ is the directed graph
 with a vertex for each \(\rho\in G^*\), and an arrow \(a_i^\rho\)
 from \(\rho\rho_i\) to \(\rho\) for each \(\rho\in G^*\) and \(1\leq
 i \leq n\).  We say the
 arrow \(a_i^\rho\) is \emph{labeled} \(i\).
 \end{definition}
 
 The McKay quiver has $r$ vertices and $nr$ arrows, and can be encoded
 in an $(r+n) \times nr$ matrix $C$ as follows.  Let $\{
 \mathbf{e}_{\rho}: \rho \in G^{*} \} \cup \{\mathbf{e}_i : 1 \leq i
 \leq n \}$ be the standard basis of $\mathbb Z^{r+n}$, and let $\{
 \mathbf{e}^{\rho}_i : \rho \in G^{*}, 1 \leq i \leq n\}$ denote the
 standard basis of $\Z^{nr}$.  Order the latter basis globally into $r$
 blocks, one for each \(\rho\in G^*\) beginning with the trivial
 representation \(\rho_0\).  Within each block the elements are listed
 \(\mathbf{e}^{\rho}_1,\dots, \mathbf{e}^{\rho}_n\).  Let $C$ be the
 $(r+n) \times nr$ matrix with column
 $\mathbf{e}_{\rho}-\mathbf{e}_{\rho \rho_i}+\mathbf{e}_i$
 corresponding to $\mathbf{e}_i^{\rho}$. Note that the top $r \times
 (nr)$ submatrix $B$ with column $\mathbf{e}_{\rho} - \mathbf{e}_{\rho
   \rho_i}$ corresponding to $\mathbf{e}_i^{\rho}$ is the vertex-edge
 incidence matrix of the McKay quiver.
 
 A \emph{representation of the McKay quiver} of dimension vector
 \((1,\dots ,1)\in \N^r\) is the assignment of a one-dimensional
 \(\K\)-vector space \(R_\rho\) to each vertex \(\rho\), and a linear
 map \(R_{\rho\rho_i}\rightarrow R_{\rho}\) to each arrow
 $a_i^{\rho}$.  Fix a basis for each \(R_{\rho}\) and write
 \(b_i^{\rho} \in \K\) for the entry of the \(1 \times 1\) matrix of
 the linear map $R_{\rho\rho_i}\rightarrow R_{\rho}$.  We occasionally
 use $b_i^{\rho}$ to refer to the linear map itself.  Since there are
 \(nr\) arrows in the quiver, representations define points
 \((b_i^\rho)\in\mathbb{A}_{\K}^{nr}\).  We write \(\K[z_i^\rho :
 \rho\in G^*, 1\leq i \leq n]\) for the coordinate ring of
 \(\mathbb{A}_{\K}^{nr}\). We consider only points \((b_i^\rho)\) of
 the scheme \(Z\) defined by the ideal \[ I = \langle z_j^{\rho\rho_i}
 z_i^{\rho} - z_i^{\rho\rho_j} z_j^{\rho} : \rho \in G^*, 1 \leq i,j
 \leq n\rangle.  \] Thus, we consider only representations
 \((b_i^\rho)\in\mathbb{A}_{\K}^{nr}\) satisfying the relations
\begin{equation} \label{e:quiverrelations}
b_j^{\rho\rho_i} b_i^{\rho} = b_i^{\rho\rho_j} b_j^{\rho} \text{ for
 }\rho\in G^* \text{ and }1\leq i,j\leq n.
\end{equation}
These relations arise naturally when quiver representations are
translated into the equivalent language of \(G\)-constellations
(see Remark~\ref{remark:ARS}).

We now summarize the Geometric Invariant Theory (GIT) construction of
the moduli spaces of $\theta$-stable McKay quiver representations (see
\cite[\S 2,\S 4]{CMT1} for more details).  The
algebraic torus $\K^*$ acts on each $R_{\rho}$, so \((\K^*)^r\) acts
diagonally on the vector space \(\oplus_{\rho\in G^*} R_\rho\) by
change of basis.  Hence \(t =(t_\rho)\in (\K^*)^r\) acts on
\(b_i^\rho\in \Hom(R_{\rho\rho_i},R_{\rho}) = R_{\rho\rho_i}^*\otimes
R_{\rho}\) as
\begin{equation} \label{eqn:rtorusaction} t \cdot
 b_i^\rho = t_{\rho\rho_i}^{-1} t_{\rho} b_i^\rho.  
\end{equation}
The diagonal scalar subgroup acts trivially, leaving a faithful action
of the $(r-1)$-dimensional algebraic torus \( T_B:= \Hom(\Z B,\K^*)\)
on \(\mathbb{A}^{nr}_{\K}\) whose character lattice \(\Z B\subset
\Z^r\) is generated by the columns of the matrix \(B\). This action
induces a $\Z B$-grading of the coordinate
  ring of \(\mathbb{A}_{\K}^{nr}\) by setting $\deg(z_i^\rho)=
  \mathbf{e}_{\rho} - \mathbf{e}_{\rho \rho_i}$.  The ideal
\(I\) defining \(Z\) is homogeneous, so $\K[Z]$ is $\Z B$-graded and,
for $\textbf{b}\in \Z B$, we write $\K[Z]_{j\textbf{b}}$ for the
  $j\mathbf{b}$-graded piece of $\K[Z]$.
  Then the categorical GIT quotient of $Z$ by the action of
  $T_B$ linearized by $\textbf{b}$ is the scheme
 \[
 Z\git_{\!\textbf{b}}T_B := \Proj \textstyle{\bigoplus_{j\geq 0}}
 \K[Z]_{j\textbf{b}}.
 \]
 More generally, the quotient linearized by an element
 \(\theta\in \Z B\otimes \Q\) in the $\Q$-vector space generated by
 the columns of $B$ is defined to be the GIT quotient linearized by
 any multiple for which \(j\theta\in \Z B\).  A parameter \(\theta\in
 \Z B\otimes \Q\) is \emph{generic} if every point of \(Z\) that is
 \(\theta\)-semistable (in the sense of GIT) is in fact
 \(\theta\)-stable, in which case \(Z\git_{\!\theta}T_B\) is a
 geometric quotient. The subset of generic parameters decomposes into
 finitely many open \emph{chambers}, where \(Z\git_{\!\theta}T_B\)
 remains unchanged as $\theta$ varies in a chamber, though its
 polarizing line bundle varies.

 \begin{definition} \label{defn:Theta} The \emph{GIT parameter space}
   is the $\Q$-vector space
 \[
 \Theta := \Z B\otimes \Q = \big{\{}(\theta_{\rho}) \in \Q^r :
 \textstyle{\sum_{\rho \in G^*} \theta_{\rho} = 0}\big{\}}. 
 \]
 For $\theta \in \Theta$, $\mathcal{M}_\theta:= Z\git_\theta T_B$ is
 the \emph{coarse moduli space of $\theta$-semistable McKay quiver
   representations} of dimension vector \((1,\dots,1)\) satisfying the
 relations~(\ref{e:quiverrelations}).  For generic $\theta$,
 $\mathcal{M}_\theta$ is the \emph{fine moduli space of
   $\theta$-stable McKay quiver representations}.
 \end{definition}
 
 The best known example of $\mathcal{M}_\theta$ is
 the \(G\)-Hilbert scheme, denoted \(\ghilb\). This parameterizes ideals
 \(J\subseteq S=\K[x_1,\dots,x_n] \) defining \(G\)-invariant
 subschemes \(Z(J)\subseteq \mathbb{A}^{n}_{\K}\) whose coordinate
 rings \(S/J\) are isomorphic to the group ring $\K G$ as \(\K G\)-modules.
 Ito--Nakajima~\cite[\S3]{ItoNakajima} observed that there is a unique
 chamber in \(\Theta\) containing parameters \(\{\theta \in \Theta \st
 \theta_\rho > 0 \text{ for } \rho\neq \rho_0\}\) such that
 \(\mathcal{M}_\theta \cong \ghilb\).

 To state the main result of Craw--Maclagan--Thomas~\cite{CMT1}, let
 $\mathbb N C\subset \Z^{r+n}$ denote the subsemigroup generated by
 the columns of the matrix $C$ and let \(P\subseteq \Q^{r+n}\) be the
 cone generated by the column vectors
 $\{\mathbf{e}_{\rho}-\mathbf{e}_{\rho \rho_i}+\mathbf{e}_i :\rho\in
 G^*, 1\leq i\leq n\}$ of \(C\).  Also, let $\pi\colon
 \Q^{r+n}\rightarrow \Q^r$ and $\pi_n\colon \Q^{r+n}\rightarrow
 \Q^n\cong \ker_{\Z}(\pi)\otimes_{\Z} \Q$ be the projections onto the
 first $r$ and last $n$ coordinates respectively.

 \begin{theorem}[Craw--Maclagan--Thomas~\cite{CMT1}]
 \label{thm:CMT1}
   The not-necessarily-normal toric variety $V = \Spec \K[\mathbb N C]$ is a $T_B$-invariant irreducible
   component of the scheme $Z\subset \mathbb{A}^{nr}_{\K}$.   In addition:
 \begin{enumerate}
 \item For $\theta \in \Theta$, the GIT quotient $Y_\theta :=
   V\git_\theta T_B$ is a not-necessarily-normal toric variety that
   admits a projective birational morphism $\tau_\theta\colon Y_\theta
   \rightarrow \mathbb{A}^{n}_{\K}/G$ obtained by variation of GIT
   quotient.
 \item For generic $\theta\in \Theta$, the variety $Y_\theta$ is the
   unique irreducible component of $\mathcal M_{\theta}$ containing
   the $T_B$-orbit closures of the points of $Z\cap (\K^*)^{nr}$.
 \item The toric fan of $Y_\theta$ is the inner normal fan of the
   polyhedron $P_\theta$ obtained as the convex hull of the set \(\pi_n(P\cap
 \pi^{-1}(\theta))\subset \ker_{\Z}(\pi)\otimes_{\Z} \Q\).
 \end{enumerate}
  \end{theorem}

 \begin{definition}
 For generic $\theta\in \Theta$, \(Y_{\theta}\) is called the
   \emph{coherent component} of \(\mathcal M_{\theta}\).
 \end{definition}
  
 In the special case where \(\mathcal{M}_\theta \cong \ghilb\), we
 established \cite[Corollary~1.2]{CMT1} that the coherent component
 $Y_\theta$ is isomorphic to the original version of the $G$-Hilbert
 scheme $\hilbg$ introduced by Nakamura~\cite{Nakamura}.

 \subsection{$G$-constellations}
 We recall the notion of $G$-constellation and review some well-known
 results from representation theory for which we could not find a
 suitable reference.

 Let $S:= \K[x_1,\dots, x_n]$. The group $G$ acts on $S$ by $g \cdot
 x_i=\rho_i(g^{-1})x_i$.  We now recall the {\em skew group algebra}
 $S \rtimes G$.  As an $S$-module, the skew group algebra is the free
 $S$-module with basis $G$.  The ring structure is given by setting
 $(sg) \cdot(s'g')=s(g \cdot s') gg'$ for $s,s'\in S$ and $g,g'\in G$.
 Recall that an $S$-module $M$ is $G$-equivariant if it has a
 $G$-action such that $g \cdot (sm)=(g \cdot s) (g \cdot m)$ for $g
 \in G, s \in S$ and $m \in M$.  An $S$-module is $G$-equivariant if
 and only if it is a left $S \rtimes G$-module.

 \begin{definition} \label{defn:gcon} A \emph{$G$-constellation} is a
  $G$-equivariant $S$-module that is isomorphic as a $\K G$-module to
  $\K G$.
 \end{definition}

 In order to apply Gr\"obner basis theory we reinterpret
 $G$-constellations as graded modules. We give $S$ a $G^*$-grading by
 $\deg(x_i)=\rho_i$ for $1 \leq i \leq n$.  This grading comes from
 the inclusion of $G$ into the $n$-dimensional torus acting on
 $\mathbb A^n_{\K}$, and gives a map $\deg\colon \Z^n \rightarrow G^*$.

\begin{definition}
Define a $G^\ast$-graded $\K$-algebra as follows.  As
a free $S$-module, $A= \oplus_{\rho \in G^*} S {\bf e}_\rho$ is the
free $S$-module of rank $r$ with basis $\{ \mathbf{e}_{\rho} : \rho \in
G^\ast \}$.  Extend the $G^*$-grading of $S$ to a $G^*$-grading on $A$
by defining $\deg(\mathbf{e}_{\rho})=\rho$.  The multiplication on $A$
is then determined by
$$\mathbf{e}_{\rho'} \cdot x^{\mathbf{u}} \mathbf{e}_{\rho} = 
\left \{ 
\begin{array}{ll}
x^{\mathbf{u}}\mathbf{e}_{\rho} & \text{ if } \deg(x^{\mathbf{u}}\mathbf{e}_{\rho})=\deg(\mathbf{e}_{\rho'}) \\
0 & {\text{ otherwise}} 
\end{array}
\right. 
$$
together with the $S$-module structure.
\end{definition}

\begin{remark}
  The algebra $A$ is the path algebra of the McKay quiver modulo the
  ideal of relations corresponding to (\ref{e:quiverrelations}).  See
  \cite[III.1]{ARS} for the definition of the path algebra.  This
  description requires the assumption that $G$ is abelian.  It is
  well-known that for finite $G$ in $\GL(n,\K)$ (see for example
  \cite[Chapter 10]{Yoshino}) the algebra $A$ is Morita equivalent to the
  skew group algebra.  In the abelian case these algebras are actually
  isomorphic.
\end{remark}

\begin{proposition} \label{p:Aisskewgroup} \begin{enumerate}
 \item An $S$-module is a left $A$-module if
  and only if it is $G^*$-graded.  
 \item An $S$-module homomorphism between
  left $A$-modules is a left $A$-module homomorphism if and only if it
  preserves the $G^*$-grading.
 \item The algebra $A$ is
  isomorphic to $S \rtimes G$. 
 \end{enumerate}
\end{proposition}

\begin{proof}
  Let $M$ be a $G^*$-graded $S$-module.  Define a left $A$-module
  structure on $M$ by setting $\mathbf{e}_{\rho} \cdot m = m_{\rho}$
  for $m \in M$, where $m_{\rho}$ is the piece of $m$ in degree
  $\rho$.  To check that this gives an $A$-module structure, it
  suffices to show that
  $(\mathbf{e}_{\rho'}(x^{\mathbf{u}}\mathbf{e}_{\rho})) \cdot m =
  \mathbf{e}_{\rho'} \cdot (x^{\mathbf{u}} \mathbf{e}_{\rho} \cdot
  m)$.  The expression on the left-hand side is
  $x^{\mathbf{u}}m_{\rho}$ if
  $\deg(\mathbf{e}_{\rho'})=\deg(x^\mathbf{u}\mathbf{e}_{\rho})$, and
  zero otherwise.  Since $x^\mathbf{u}\mathbf{e}_{\rho} \cdot
  m=x^\mathbf{u} m_{\rho}$, this equals the right-hand side.
 
 Conversely, let $M$ be an $A$-module.  Define
  $M_{\rho}=\mathbf{e}_{\rho} M$.  We claim that $M=\oplus_{\rho \in
  G^*} M_{\rho}$ as an abelian group, and this decomposition is
  compatible with multiplication by elements of $S$.  Indeed, if $m
  \in M_{\rho} \cap M_{\rho'}$, then $m=\mathbf{e}_{\rho}
  m_1=\mathbf{e}_{\rho'} m_2$ for some $m_1, m_2 \in M$.  But then
  $\mathbf{e}_{\rho} m = \mathbf{e}_{\rho}^2m_1=\mathbf{e}_{\rho}
  m_1=m$, and
  $\mathbf{e}_{\rho}m=\mathbf{e}_{\rho}\mathbf{e}_{\rho'}m_2=0$, so
  $m=0$.  Since $e = \sum_{\rho \in G^*} \mathbf{e}_{\rho}$ is the
  multiplicative identity of $A$, we have $m=e\cdot m=\sum_{\rho \in
  G^*} m_{\rho}$.  This gives the decomposition as abelian groups. If
  $\deg(x^\mathbf{u})=\rho'$ then $\mathbf{e}_{\rho' \rho} \cdot
  x^\mathbf{u} \mathbf{e}_{\rho}m_{\rho} = x^\mathbf{u}m_{\rho}$.
  This gives $x^\mathbf{u} m_{\rho} \in M_{\rho' \rho}$, so (1) holds.

An $A$-module homomorphism $\phi$ satisfies $\phi(\mathbf{e}_{\rho} m)
= \mathbf{e}_{\rho} \phi(m)$, so is exactly a degree-zero $G^*$-graded
$S$-module homomorphism.  This gives (2).

 The $\K$-linear map $\phi \colon \oplus_{\rho \in G^*}
  \K\mathbf{e}_{\rho} \rightarrow \oplus_{g \in G} \K g$ given by
  $\phi(\mathbf{e}_{\rho})=1/r \sum_{g \in G} \rho(g)g$ is an
  isomorphism since the character table is an invertible matrix for an
  abelian group.  This extends to an isomorphism of $S$-modules $\phi
  \colon A \rightarrow S \rtimes G$.  It remains to check that
  $\phi(\mathbf{e}_{\rho})\phi(x^{\mathbf{u}}\mathbf{e}_{\rho'}) =
  \phi(\mathbf{e}_{\rho} \cdot x^{\mathbf{u}} \mathbf{e}_{\rho'})$.
  Indeed, the left-hand side is
  $$
  \frac{1}{r^2} x^{\mathbf{u}} \sum_{h \in G} \sum_{g,g' \in G,\\
    gg'=h} (\rho \rho'')(g) \rho'(g')h,$$
  where
  $\deg(x^{\mathbf{u}})={\rho''}^{-1}$, so $g\cdot
  x^{\mathbf{u}}=\rho''(g)x^{\mathbf{u}}$.  By the orthogonality
  relations of the character table, this is $x^{\mathbf{u}}/r \sum_{g
    \in G} \rho'(g) g$ when $\rho' = \rho\rho''$, and zero otherwise.
  This proves the final statement.
\end{proof}

 \begin{corollary}
 \label{c:gconisgraded}
 An $S$-module $F$ is a $G$-constellation if and only if $F$ is $G^*$-graded with
 Hilbert function $\dim_{\K} F_\rho = 1$ for each degree $\rho\in
 G^*$.
 \end{corollary}

 \begin{proof} 
   A $G$-constellation is an $S \rtimes G$-module, and hence is
   $G^*$-graded by Proposition~\ref{p:Aisskewgroup}. It remains to
   show that a $G^*$-graded module is isomorphic to $\K G=\oplus_{g
     \in G} \K g$ if and only if it has Hilbert function one in each
   degree.  This follows from the $\K$-linear isomorphism $\phi$ from
   the proof of Proposition~\ref{p:Aisskewgroup}, and the fact that,
   if $F$ is a $G^*$-graded module with Hilbert function one in each
   degree, there is a surjection from $\oplus_\rho S {\bf e}_\rho$
   where the image of each ${\bf e}_\rho$ is nonzero.
\end{proof}


\section{Gr\"obner interpretation of $G$-constellations}
\label{sec:grobner}
In this section we first review some well-known facts from the theory of
Gr\"{o}bner basis for modules. We then canonically associate a submodule of $A$
to every $G$-constellation, and establish the key Gr\"{o}bner
result by exhibiting a Gr\"{o}bner basis for this module.

\subsection{Preliminary Gr\"obner facts}
\label{sec:prelimgrobner}
We start by summarizing the relevant facts about Gr\"obner bases (see
Cox--Little--O'Shea~\cite{CLO} and Eisenbud~\cite[Chapter
15]{Eisenbud}).

Let $M$ be a submodule of the free module $S^r$ for some $r \in
\mathbb N$.  An element $f \in M$ can be written as $f=\sum
c_{\mathbf{u},i} x^{\mathbf{u}}\mathbf{e}_{i}$, where the sum is over
$\mathbf{u} \in \mathbb N^n$ and $1 \leq i \leq r$, and all but
finitely many $c_{\mathbf{u},i}$ are zero. If
$\mathbf{v}=(v_1,\dots,v_r) \in \mathbb Q^r$, and $\mathbf{w} \in
\mathbb Q_{\geq 0}^n$, then the initial term
$\inn_{(\mathbf{v},\mathbf{w})}(f)=\sum c_{\mathbf{u},i}
x^{\mathbf{u}} \mathbf{e}_{i}$, where the sum is over pairs
$(\mathbf{u},i)$ with $\mathbf{w} \cdot \mathbf{u} + v_i
\geq\mathbf{w} \cdot \mathbf{u}' + v_j$ for any other pair
$(\mathbf{u}',j)$ with $c_{\mathbf{u}',j} \neq 0$.  The initial module
of $M$ is the $S$-module
$$\inn_{(\mathbf{v},\mathbf{w})}(M)=\langle
\inn_{(\mathbf{v},\mathbf{w})}(f) : f \in M \rangle.$$ 

If $M$ is homogeneous in some grading of $S^r$, then the Hilbert
function of $S^r/M$ equals that of
$S^r/\inn_{(\mathbf{v},\mathbf{w})}(M)$.  When
$\inn_{(\mathbf{v},\mathbf{w})}(M)$ is generated by monomials
$x^{\mathbf{u}}\mathbf{e}_i$ then a basis for
$S^r/\inn_{(\mathbf{v},\mathbf{w})}(M)$ consists of those monomials
not lying in $\inn_{(\mathbf{v},\mathbf{w})}(M)$, which we call the
{\em standard monomials} of $\inn_{(\mathbf{v},\mathbf{w})}(M)$.

It is important to note that the initial module cannot usually be
computed by taking the initial terms of the module generators.  A {\em
Gr\"obner basis} with respect to the weight vector $(\mathbf{v},
\mathbf{w})$ is a set of generators $\{ m_1, \dots ,m_s \}$ for $M$
such that $\inn_{(\mathbf{v},\mathbf{w})}(M)=\langle
\inn_{(\mathbf{v},\mathbf{w})}(m_1), \dots,
\inn_{(\mathbf{v},\mathbf{w})}(m_s) \rangle$.

Gr\"obner bases are usually defined by giving a {\em term order},
which is a total order on the monomials $x^{\mathbf{u}}\mathbf{e}_i$
in $S^r$ satisfying $\mathbf{e}_i \prec x^{\mathbf{u}} \mathbf{e}_i$
for all $\mathbf{u} \neq 0$, and if $x^{\mathbf{u}} \mathbf{e}_i \prec
x^{\mathbf{u}'} \mathbf{e}_j$ then
$x^{\mathbf{u}+\mathbf{u}''}\mathbf{e}_i \prec
x^{\mathbf{u}'+\mathbf{u}''} \mathbf{e}_j$.  A weight vector
$(\mathbf{v},\mathbf{w})$ gives a partial order, called a {\em weight
  order}, by setting $x^{\mathbf{u}} \mathbf{e}_i \prec
x^{\mathbf{u}'} \mathbf{e}_j$ if $\mathbf{w} \cdot {\mathbf{u}} + v_i
< \mathbf{w} \cdot \mathbf{u}' + v_j$.  This partial order can be
refined to a term order by breaking ties with a fixed term order.  We
use only the {\em term over position} lexicographic
order~\cite[Chapter 5, Definition 2.4]{CLO2}.  We will use the
following proposition (see Sturmfels~\cite[Corollary 1.9]{GBCP} for a
proof for ideals in a polynomial ring; the extension to modules is
straightforward).

\begin{proposition} \label{p:weightvect}
Let $(\mathbf{v},\mathbf{w})$ be a weight vector and
$\prec_{(\mathbf{v},\mathbf{w})}$ be a term order that refines the
weight order.  If $\{ m_1, \dots m_s\}$ is a Gr\"obner basis for a
module $M$ with respect to $\prec_{(\mathbf{v},\mathbf{w})}$, then
$\{m_1,\dots, m_s \}$ is also a Gr\"obner basis for the weight order
given by $(\mathbf{v},\mathbf{w})$.
\end{proposition}

A criterion for a subset $\{m_1,\dots,m_s\}$ of $M$ to be a Gr\"obner
basis is given by the conditions of {\em Buchberger's algorithm}.  The
key idea in this algorithm is that of an {\em S-pair}: if
$\inn_{\prec}(m_i)=c_ix^{\mathbf{u}_i}\mathbf{e}_k$ and
$\inn_{\prec}(m_j)=c_jx^{\mathbf{u}_j}\mathbf{e}_k$ involve the same
basis element $\mathbf{e}_k$ of $S^r$, then $S(m_i,m_j):=m_{ji}
m_i-m_{ij} m_j$, where $m_{ji}=c_j
x^{\mathbf{u}_j}/\gcd(x^{\mathbf{u}_i},x^{\mathbf{u}_j})$ and
$m_{ij}=c_i x^{\mathbf{u}_i}/\gcd(x^{\mathbf{u}_i},x^{\mathbf{u}_j})$.
The set $\{m_1,\dots,m_s\}$ is a Gr\"obner basis if every $S$-pair can
be written as $S(m_i,m_j)=\sum_l h_l m_l$, where $h_l \in S$ and
$\inn_{\prec}(h_lm_l) \preceq \inn_{\prec}(S(m_i,m_j))$.  In this case
we say that the $S$-pair {\em reduces to zero}.  In general, if
$f=\sum_l h_l m_l +g$, where $\inn_{\prec}(h_lm_l) \preceq
\inn_{\prec}(f)$ then we say that $f$ reduces to $g$ modulo
$\{m_1,\dots, m_s \}$.

 \subsection{The key Gr\"{o}bner basis result}
 We now use Gr\"{o}bner basis techniques to write down an explicit map
 that associates to each quiver representation $(b_i^{\rho}) \in Z$, a
 $G$-constellation with a chosen presentation.  This map arises
 naturally via an isomorphism of categories.  First, we introduce the
 categories.

\begin{definition} \label{d:categoryR}
Let $\mathcal R$ be the category whose objects are points $(b_i^{\rho})
\in Z$.  A morphism $h \colon (b_i^{\rho}) \rightarrow
({b'}_i^{\rho})$ consists of a scalar $h_{\rho} \in \K$ for every $\rho
\in G^*$, satisfying 
$b_i^{\rho} h_{\rho} = h_{\rho \rho_i}
{b'}_i^{\rho}$.
\end{definition}

\begin{remark}
The category $\mathcal R$ is obtained from the category of McKay
quiver representations of dimension vector $(1, \dots, 1)$ satisfying
the given relations (see \cite[\S III]{ARS}) by choosing a basis for
each $R_{\rho}$.  In particular,  the commutative diagram
$$\xymatrix{R_{\rho \rho_i} \ar[r]^{b_i^{\rho}} \ar[d]_{h_{\rho \rho_i}} & R_{\rho}
\ar[d]^{h_{\rho}} \\ R'_{\rho \rho_i} \ar[r]^{{b'}_i^{\rho}} & R'_{\rho} \\ }
$$
leads naturally to the conditions $b_i^{\rho} h_{\rho} = h_{\rho
  \rho_i} {b'}_i^{\rho}$.
\end{remark}

\begin{definition}
  Let $\mathcal C$ be the category whose objects are cyclic $A$-modules of the
  form $A/M$ satisfying $\dim_{\K}(A/M)_{\rho}=1$ for all $\rho \in
  G^*$, with $\mathbf{e}_{\rho} \not \in M$ for our chosen $S$-module
  basis $\mathbf{e}_{\rho}$ of $A$.   The morphisms of
  $\mathcal{C}$ are
  $A$-module homomorphisms.
\end{definition}

\begin{remark}
Corollary~\ref{c:gconisgraded} shows that the objects of $\mathcal C$
are $G$-constellations with a chosen presentation.
\end{remark}

We now construct an explicit isomorphism of categories between
$\mathcal C$ and $\mathcal R$.  We emphasize that it is
highly unusual in Gr\"obner theory for the minimal generating set of a
module to be a Gr\"obner basis, as occurs in this proposition.

  \begin{proposition} \label{p:mggcons} There is a contravariant
functor $\Psi: \mathcal R \rightarrow \mathcal C$ taking the McKay
quiver representation $b = (b_i^\rho)\in Z$ to the $G$-constellation
$A/M_b$ for the left $A$-ideal
 \begin{equation}
 \label{eqn:Mb}
 M_b:= \langle x_i\mathbf{e}_{\rho }-b_i^{\rho}
 \mathbf{e}_{\rho \rho_i} : 1 \leq i \leq n, \rho \in G^* \rangle
 \end{equation}
 that is an isomorphism of categories.  Moreover, the given generators
 for $M_{b}$ form a Gr\"obner basis with respect to any monomial order
 with $x_i \mathbf{e}_{\rho } \succ \mathbf{e}_{\rho \rho_i}$.
\end{proposition}

 \begin{proof} We begin by showing that the given generators for
   $M_{b}$ form a Gr\"obner basis for any term order with $x_i
   \mathbf{e}_{\rho } \succ \mathbf{e}_{\rho \rho_i}$.  Indeed, the
   only relevant $S$-pairs are between terms of the form
   $x_i\mathbf{e}_{\rho }-b_i^{\rho }\mathbf{e}_{\rho \rho_i}$ and
   $x_j \mathbf{e}_{\rho}- b_j^{\rho \rho_j} \mathbf{e}_{\rho
   \rho_j}$.  This $S$-pair is then $x_j b_i^{\rho} \mathbf{e}_{\rho
   \rho_i}-x_i b_j^{\rho }\mathbf{e}_{\rho \rho_j}$.  These terms are
   symmetric in $i,j$, so we may assume that the first term is the
   leading term in our term order, and that $b_i^{\rho} \neq 0$.  The
   polynomial now reduces using the binomial $x_j \mathbf{e}_{\rho
   \rho_i} - b_j^{\rho \rho_i} \mathbf{e}_{\rho \rho_i \rho_j}$ to
   $b_j^{\rho \rho_i} b_i^{\rho} \mathbf{e}_{\rho \rho_i \rho_j} - x_i
   b_j^{\rho} \mathbf{e}_{\rho \rho_j}$.  Our assumption on the term
   order now implies that if $b_j^{\rho} \neq 0$, the second term is
   the leading term, so the binomial reduces using $x_i
   \mathbf{e}_{\rho \rho_j} - b_i^{\rho \rho_j} \mathbf{e}_{\rho \rho_i \rho_j}$ to
   $(b_i^{\rho \rho_j} b_j^{\rho} - b_j^{\rho \rho_i} b_i^{\rho})
   \mathbf{e}_{\rho}$.  This is zero since $(b_i^{\rho}) \in Z$.  If
   $b_j^{\rho} =0$, then since $b_i^{\rho} \neq 0$ and $(b_i^{\rho})
   \in Z$ we must have $b_j^{\rho \rho_i}=0$, so the intermediate
   binomial was already zero.  In both cases the $S$-pair reduces to
   zero, so the generators for $M_b$ form a Gr\"obner basis.

   We next show that $A/M_{b}$ is an object of $\mathcal C$.  Since
   $\deg(x_i \mathbf{e}_{\rho})=\deg(\mathbf{e}_{\rho \rho_i})$, the
   submodule $M_{b}$ is homogeneous in the $G^*$-grading, so $A/M_{b}$
   is also graded by $G^*$.  Since the given generators for $M_{b}$
   form a Gr\"obner basis as above, $M_0=\langle x_i \mathbf{e}_{\rho}
   : 1 \leq i \leq n, \rho \in G^* \rangle$ is an initial module of
   $M_{b}$, and thus the Hilbert function of $A/M_0$ equals that of
   $A/M_{b}$.  The Hilbert function of the quotient by a monomial
   module is the number of standard monomials of the module in the
   given degree.  Since the only standard monomials of $M_0$ are the
   units $\mathbf{e}_{\rho}$, we conclude that the Hilbert function of
   $A/M_b$ is one in every degree.

We now construct $\Psi$ and its inverse.  Define $\Psi$ on objects as
above.  A morphism $h \colon (b_i^{\rho}) \rightarrow ({b'}_i^{\rho})$
is a collection of $h_{\rho} \in \K$ for $\rho \in G^*$.  Define an
$S$-module homomorphism $\psi(h):A \rightarrow A/M_{b}$ by
$\psi(h)(\mathbf{e}_{\rho})=h_{\rho} \mathbf{e}_{\rho}$.  Since
$\psi(h)(x_i\mathbf{e}_{\rho }-{b'}_i^{\rho} \mathbf{e}_{\rho \rho_i}) =
h_{\rho } (x_i\mathbf{e}_{\rho}-b_i^{\rho}
\mathbf{e}_{\rho \rho_i}) \in M_{b}$, $\psi(h)$ defines an $S$-module
homomorphism $\Psi(h) \colon A/M_{b'} \rightarrow A/M_{b}$.  It is now
straightforward to check that $\Psi$ is a functor from $\mathcal R$ to
$\mathcal C$.  To construct the inverse functor $\Phi$, let $A/M$ be
an object of $\mathcal C$.  Since $\dim_{\K}(A/M)_{\rho}=1$ for all
$\rho \in G^*$, $\deg(x_i\mathbf{e}_{\rho})=\deg(\mathbf{e}_{\rho \rho_i})$, and $\mathbf{e}_{\rho} \not \in M$,
there is a unique $b_i^{\rho}\in \K$ with $x_i \mathbf{e}_{\rho} -
b_i^{\rho} \mathbf{e}_{\rho \rho_i} \in M$ for each $i$ and $\rho$.  We then
have $x_ix_j \mathbf{e}_{\rho} -b_i^{\rho \rho_j}
b_j^{\rho} \mathbf{e}_{\rho \rho_i \rho_j} \in M$, and hence $(b_i^{\rho
\rho_j}b_j^{\rho} -b_j^{\rho \rho_i}b_i^{\rho}) \mathbf{e}_{\rho \rho_i \rho_j} \in
M$.  Since $\mathbf{e}_{\rho} \not \in M$, we conclude that $b_i^{\rho
\rho_j}b_j^{\rho} -b_j^{\rho \rho_i}b_i^{\rho}=0$, and so
$(b_i^{\rho}) \in Z$.  We thus set $\Phi(A/M)=(b_i^{\rho})$.  Given a
morphism $h:A/M \rightarrow A/M'$, lift to $\tilde{h} : A \rightarrow
A/M'$, and write $\tilde{h}(\mathbf{e}_{\rho})=\lambda_{\rho}
x^{\mathbf{u}} \mathbf{e}_{\rho'}$, for some $u, \rho'$ satisfying
$\deg(x^{\mathbf{u}} \mathbf{e}_{\rho'})=\deg(\mathbf{e}_{\rho})$.
Since $\dim_{\K}(A/M')_{\rho}=1$, and $\mathbf{e}_{\rho} \not \in
M'$, there is a unique $\mu \in \K$ with
$x^{\mathbf{u}}\mathbf{e}_{\rho \rho'}-\mu \mathbf{e}_{\rho} \in M'$,
and thus $h(\mathbf{e}_{\rho})=\lambda_{\rho} \mu \mathbf{e}_{\rho}$.
The scalar $\lambda_{\rho} \mu$ is independent of the choice of $u$
and $\rho'$ since $\mathbf{e}_{\rho} \not \in M'$.  We define
$\Phi(h)$ to be the morphism $\Phi(h) \colon ({b'}_i^{\rho}) \rightarrow
(b_i^{\rho})$ with $\Phi(h)_{\rho}=\lambda_{\rho} \mu$.  The fact
$\tilde{h}(M)=0$ implies that $\Phi(h)$ is a morphism in $\mathcal R$.
It follows that $\Phi$ is a functor from $\mathcal C$ to $\mathcal R$,
and $\Phi=\Psi^{-1}$.
\end{proof}

\begin{remark}
 \label{remark:ARS}
 \begin{enumerate}
 \item The \(A\)-submodules \(M_b\subseteq A\) may be regarded
 as left ideals in the skew group ring \(A \cong S\rtimes G\).  
 \item The relations~(\ref{e:quiverrelations}) correspond via $\Psi$  to the
 commutativity of $x_i$ and $x_j$ in the $S$-module structure on
 $A/M_b$.
\item The translation from quiver representations to modules over an
  algebra is a special case of a result for representations of
  an arbitrary finite quiver with relations (see
  \cite[III,Proposition~1.7]{ARS}). Since $\mathcal C$ and
  $\mathcal R$ involve choices of bases, we obtain an isomorphism
  rather than an equivalence of categories.
\item The isomorphism of categories $\Psi$ implies that
  \(\mathcal{M}_\theta\) from Definition~\ref{defn:Theta} may be
  regarded as the moduli space of $\theta$-semistable
  $G$-constellations.
 \end{enumerate}
\end{remark}

\section{Distinguished \protect\(G\protect\)-constellations via Gr\"obner bases}
 \label{sec:distinguishedGcons}
 This section explicitly describes the distinguished
 $G$-constellations that define distinguished points on
 $Y_\theta\subseteq \mathcal{M}_\theta$.  We exploit here the
 geometric interpretation of Gr\"obner bases as allowing explicit
 computations of flat degenerations coming from a one-parameter torus
 (see \cite[Chapter 15]{Eisenbud}).  This Gr\"{o}bner description
 allows us to decide whether or not a given point of
 $\mathcal{M}_\theta$ lies in $Y_\theta$.

\subsection{Distinguished \protect\(G\protect\)-constellations}
In order to apply the Gr\"{o}bner result from Proposition~\ref{p:mggcons},  we introduce an $A$-module that plays a key role in the rest of the paper.

 \begin{definition}
 \label{defn:McKaymod}
 The {\em McKay module} is the submodule, or left ideal, of $A$ given
   by $$ M_G = \langle x_i \mathbf{e}_{\rho}- \mathbf{e}_{\rho\rho_i}
   : \rho\in G^*, 1\leq i\leq n\rangle.$$ Note that $M_G = M_b$ for
   the point $b\in Z$ with $b_i^\rho = 1$ for all $\rho\in G^*$ and
   $1\leq i\leq n$.
 \end{definition}

 \begin{lemma}
   The McKay module $M_G$ is equal to the $A$-module $\langle x_i - e : 1\leq
   i\leq n\rangle$, where $e$ is the multiplicative
   identity in $A$.
 \end{lemma}
 \begin{proof}
   Write $M_G' := \langle x_i - e : 1\leq i\leq n\rangle$ and fix
   \(1\leq i\leq n\).  Since $e = \sum_{\rho\in G^*} {\bf e}_\rho$, we
   have $x_ie - e = \sum_{\rho\in G^*} x_i{\bf e}_\rho -
   \sum_{\rho'\in G^*}{\bf e}_{\rho'}$.  Relabel $\rho' = \rho\rho_i$
   and regard the second term as a sum over $\rho\in G^*$ to give
   $x_ie - e = \sum_{\rho\in G^*} (x_i{\bf e}_\rho - {\bf
     e}_{\rho\rho_i})$, hence $M_G'\subseteq M_G$.  For the opposite
   inclusion, note that for every $\rho\in G^*$ and $1\leq i\leq n$ we
   have ${\bf e}_{\rho\rho_i}\cdot x_ie = x_i{\bf e}_{\rho}$.  This
   implies that $x_i{\bf e}_{\rho} - {\bf e}_{\rho\rho_i} ={\bf
     e}_{\rho\rho_i}\cdot x_ie - {\bf e}_{\rho\rho_i} = {\bf
     e}_{\rho\rho_i}(x_i-e)$, so $M_G\subseteq M_G'$ as required.
 \end{proof}

 Recall from Section~\ref{sec:moduli} that \(P\subseteq \Q^{r+n}\) is
 the polyhedral cone generated by the vectors
 $\{\mathbf{e}_{\rho}-\mathbf{e}_{\rho \rho_i}+\mathbf{e}_i :\rho\in
 G^*, 1\leq i\leq n\}$. For ${\bf w} = (w_i) \in (\mathbb Q^n_{\geq
   0})^*$, let \(P^\vee_{\bf w} := \{\mathbf{v} \in (\Q^{r})^* :
 w_i+v_\rho - v_{\rho\rho_i}\geq 0\}\) denote the slice of the dual
 cone \(P^\vee\).  For \(\theta\in \Theta\), the vector ${\bf w}\in
 (\mathbb Q^n_{\geq 0})^*$ determines a unique distinguished point
 \([b_{\theta,{\bf w}}]\in Y_\theta\) as described in
 \cite[Section~7]{CMT1}, and hence a \emph{distinguished
   \(\theta\)-semistable \(G\)-constellation} $\Psi(b_{\theta,{\bf
     w}})$ that we denote $A/M_{\theta,{\bf w}}$.

 \begin{theorem}
 \label{thm:main}
  For $\theta \in \Theta$ and $\mathbf{w} \in (\mathbb Q^n_{\geq 0})^*$, let $\mathbf{v} \in P^{\vee}_{\mathbf{w}}$ be any vector satisfying $\theta \cdot
   \mathbf{v} \leq \theta \cdot \mathbf{v}'$ for all $\mathbf{v}' \in
   P^{\vee}_{\mathbf{w}}$. Then the following $G$-constellations coincide:
 \begin{enumerate}
 \item The distinguished \(\theta\)-semistable $G$-constellation $A/M_{\theta,{\bf w}}$;
 \item The cyclic \(A\)-module $A/M_b$, where $M_b\subset A$ is the left $A$-ideal generated by $\{x_i\mathbf{e}_{\rho }-b_i^{\rho}\mathbf{e}_{\rho \rho_i} : 1 \leq i \leq n, \rho \in G^* \}$,  and $b = (b_i^\rho)$ satisfies
 \begin{equation}
 \label{eqn:birho}
 b_i^\rho =
 \left\{\begin{array}{cl} 1 & \text{if } w_i+v_\rho - v_{\rho\rho_i} =
     0 \\ 0 & \text{if } w_i+v_\rho - v_{\rho\rho_i} > 0
 \end{array}\right.;
 \end{equation}
\item The cyclic \(A\)-module $A/\inn_{(\mathbf{v},\mathbf{w})}(M_G)$,
  where $\inn_{({\bf v},{\bf w})}(M_G)$ is the {\em initial module} of
  $M_G$ with respect to $({\bf v},{\bf w})$.
 \end{enumerate}
 \end{theorem}
\begin{proof} 
  For $\theta \in \Theta$ and $\mathbf{w} \in (\mathbb Q^n_{\geq
    0})^*$, fix $\mathbf{v} \in P^{\vee}_{\mathbf{w}}$ with $\theta
  \cdot \mathbf{v} \leq \theta \cdot \mathbf{v}'$ for $\mathbf{v}' \in
  P^{\vee}_{\mathbf{w}}$.  The coordinates of the distinguished
  \(\theta\)-semistable quiver representation \(b =b_{\theta,{\bf
      w}}\) satisfy the conditions from the second part of
  Theorem~\ref{thm:main} by
  Craw--Maclagan--Thomas~\cite[Theorem~7.2]{CMT1}.  Then the
  $G$-constellation $A/M_{\theta,{\bf w}}:= \Psi(b_{\theta,{\bf w}})$
  coincides with that from (2) by Proposition~\ref{p:mggcons}.

  To see that (2) and (3) coincide, let
  $\prec_{\mathbf{v},\mathbf{w}}$ be the term order on $A$ refining
  the weight order given by $(\mathbf{v},\mathbf{w})$, where ties are
  broken using the term over position lexicographic order.  Since
  $(\mathbf{v},\mathbf{w}) \in P^{\vee}$, we have $w_i+v_{\rho} \geq
  v_{\rho\rho_i}$, and hence
  $\inn_{\prec_{\mathbf{v},\mathbf{w}}}(x_i\mathbf{e}_{\rho} -
  \mathbf{e}_{\rho\rho_i} ) = x_i \mathbf{e}_{\rho}$.
  Proposition~\ref{p:mggcons} states that the generators $\{ x_i
  \mathbf{e}_{\rho} - \mathbf{e}_{\rho\rho_i} \}$ for $M_G$ are a
  Gr\"obner basis for the term order $\prec_{({\bf v},{\bf w})}$, and
  by Proposition~\ref{p:weightvect} they also form a Gr\"obner basis
  for the weight order given by $({\bf v},{\bf w})$. This means that
  $\inn_{(\mathbf{v},\mathbf{w})}(M_G)=\langle x_i \mathbf{e}_{\rho
  }-\mathbf{e}_{\rho\rho_i } : w_i+v_{\rho}=v_{\rho \rho_i} \rangle +
  \langle x_i \mathbf{e}_{\rho} : w_i + v_{\rho} > v_{\rho \rho_i}
  \rangle = M_b$ for $b = (b_i^\rho)$ satisfying (\ref{eqn:birho}).
  This completes the proof.
 \end{proof}

 Theorem~\ref{thm:main} gives an algorithm to compute the
 distinguished $\theta$-semistable $G$-constellation $A/M_{\theta,{\bf
     w}}$.  This is the $G$-constellation
 analogue of \cite[Algorithm~7.6]{CMT1}.

\begin{algorithm}
  \label{alg:distinguished} To compute the distinguished
  $G$-constellation $A/M_{\theta,{\bf w}}$.
 \end{algorithm}

\noindent{\bf Input:} $(\theta, {\bf w}) \in \Theta \times (\mathbb
    Q^n_{\geq 0})^{\ast}$ and the McKay module $M_G$.

 \medskip

\begin{enumerate}
\item Set $M_G = \langle x_i \mathbf{e}_{\rho} - \mathbf{e}_{\rho
    \rho_i} : 1 \leq i \leq n, \rho \in G^{\ast} \rangle$.
\item Compute an optimal solution ${\bf v}$ of the linear program
$$ \textup{minimize} \{ \theta \cdot {\bf v}' \,:\, {\bf v}' \in
P^\vee_{\bf w} \}.$$
\item Compute $\inn_{({\bf v},{\bf w})}(M_G) = \langle \inn_{({\bf
      v},{\bf w})}(x_i \mathbf{e}_{\rho} - \mathbf{e}_{\rho \rho_i})
  \,:\, \rho\in G^*, 1\leq i\leq n\rangle$.  Then $A/M_{\theta, {\bf
      w}}$ has $ M_{\theta, {\bf w}} = \inn_{({\bf v},{\bf w})}(M_G).$
\end{enumerate}

\begin{proof}[Proof of Correctness]
  This is immediate from Theorem~\ref{thm:main}.
\end{proof} 

 \begin{example}
 \label{ex:hard}
 Consider the group $G \cong Z / 11 \mathbb Z$ generated by
 $\diag(\zeta, \zeta^2, \zeta^8)$.  The given three-dimensional
 representation decomposes as $\rho_1\oplus\rho_2\oplus\rho_{8}$, so
 $$ M_G = \langle x_1\mathbf{e}_{\rho_0} - \mathbf{e}_{\rho_1},
 x_2\mathbf{e}_{\rho_0} - \mathbf{e}_{\rho_2}, x_3\mathbf{e}_{\rho_0}
 - \mathbf{e}_{\rho_8}, x_1\mathbf{e}_{\rho_1} - \mathbf{e}_{\rho_2},
 x_2\mathbf{e}_{\rho_1} - \mathbf{e}_{\rho_3}, x_3\mathbf{e}_{\rho_1}
 - \mathbf{e}_{\rho_9},
 $$
 $$
 x_1\mathbf{e}_{\rho_2} - \mathbf{e}_{\rho_3}, x_2\mathbf{e}_{\rho_2} - \mathbf{e}_{\rho_4}, x_3\mathbf{e}_{\rho_2}
 - \mathbf{e}_{\rho_{10}}, x_1\mathbf{e}_{\rho_3} - \mathbf{e}_{\rho_4}, x_2\mathbf{e}_{\rho_3} -
 \mathbf{e}_{\rho_5}, x_3\mathbf{e}_{\rho_3} - \mathbf{e}_{\rho_0}, x_1\mathbf{e}_{\rho_4} - \mathbf{e}_{\rho_5},
 $$
 $$
 x_2\mathbf{e}_{\rho_4} - \mathbf{e}_{\rho_6}, x_3\mathbf{e}_{\rho_4} - \mathbf{e}_{\rho_1},
 x_1\mathbf{e}_{\rho_5} - \mathbf{e}_{\rho_6}, x_2\mathbf{e}_{\rho_5} - \mathbf{e}_{\rho_7}, x_3\mathbf{e}_{\rho_5}
 - \mathbf{e}_{\rho_2}, x_1\mathbf{e}_{\rho_6} - \mathbf{e}_{\rho_7}, x_2\mathbf{e}_{\rho_6} -
 \mathbf{e}_{\rho_8},$$
 $$
 x_3\mathbf{e}_{\rho_6} - \mathbf{e}_{\rho_3}, x_1\mathbf{e}_{\rho_7} -
 \mathbf{e}_{\rho_8}, x_2\mathbf{e}_{\rho_7} - \mathbf{e}_{\rho_9}, x_3\mathbf{e}_{\rho_7} - \mathbf{e}_{\rho_4},
 x_1\mathbf{e}_{\rho_8} - \mathbf{e}_{\rho_9}, x_2\mathbf{e}_{\rho_8} - \mathbf{e}_{\rho_{10}},
 x_3\mathbf{e}_{\rho_8} - \mathbf{e}_{\rho_5},$$
 $$
 x_1\mathbf{e}_{\rho_9} - \mathbf{e}_{\rho_{10}},
 x_2\mathbf{e}_{\rho_9} - \mathbf{e}_{\rho_0}, x_3\mathbf{e}_{\rho_9} - \mathbf{e}_{\rho_6},
 x_1\mathbf{e}_{\rho_{10}} - \mathbf{e}_{\rho_0}, x_2\mathbf{e}_{\rho_{10}} - \mathbf{e}_{\rho_1},
 x_3\mathbf{e}_{\rho_{10}} - \mathbf{e}_{\rho_7} \rangle.  $$

 We compute a distinguished \(\theta\)-stable quiver representation
 for the parameter $\theta = (1,1,1,1,-7,-9,1,1,1,8,1)$ (compare
 \cite[Example~7.7]{CMT1}).  The vector ${\bf w} = (10,7,6)$ lies in
 the relative interior of a three-dimensional cone in the fan of
 $Y_\theta$, and hence defines a torus-invariant $G$-constellation.
 The vector ${\bf v} = (-8,-10,-1, -3,6,4, -9, 0, -2, -15, -6)$ is an
 optimal solution to the linear program in Step (2) from
 Algorithm~\ref{alg:distinguished}, and Step (3) gives $$
 M_{\theta,{\bf w}} = \langle x_1\mathbf{e}_{\rho_0},
 x_2\mathbf{e}_{\rho_0} - \mathbf{e}_{\rho_2}, x_3\mathbf{e}_{\rho_0}
 - \mathbf{e}_{\rho_8}, x_1\mathbf{e}_{\rho_1}, x_2\mathbf{e}_{\rho_1}
 - \mathbf{e}_{\rho_3}, x_3\mathbf{e}_{\rho_1},
 x_1\mathbf{e}_{\rho_2}, x_2\mathbf{e}_{\rho_2} - \mathbf{e}_{\rho_4},
 $$
 $$
 x_3\mathbf{e}_{\rho_2}, x_1\mathbf{e}_{\rho_3},
 x_2\mathbf{e}_{\rho_3} - \mathbf{e}_{\rho_5}, x_3\mathbf{e}_{\rho_3},
 x_1\mathbf{e}_{\rho_4}, x_2\mathbf{e}_{\rho_4},
 x_3\mathbf{e}_{\rho_4}, x_1\mathbf{e}_{\rho_5},
 x_2\mathbf{e}_{\rho_5}, x_3\mathbf{e}_{\rho_5},
 x_1\mathbf{e}_{\rho_6}, $$
 $$
 x_2\mathbf{e}_{\rho_6} -
 \mathbf{e}_{\rho_8}, x_3\mathbf{e}_{\rho_6} - \mathbf{e}_{\rho_3},
 x_1\mathbf{e}_{\rho_7}, x_2\mathbf{e}_{\rho_7},
 x_3\mathbf{e}_{\rho_7} - \mathbf{e}_{\rho_4}, x_1\mathbf{e}_{\rho_8},
 x_2\mathbf{e}_{\rho_8}, x_3\mathbf{e}_{\rho_8} - \mathbf{e}_{\rho_5},
 $$
 $$
 x_1\mathbf{e}_{\rho_9}, x_2\mathbf{e}_{\rho_9} -
 \mathbf{e}_{\rho_0}, x_3\mathbf{e}_{\rho_9} - \mathbf{e}_{\rho_6},
 x_1\mathbf{e}_{\rho_{10}}, x_2\mathbf{e}_{\rho_{10}},
 x_3\mathbf{e}_{\rho_{10}} - \mathbf{e}_{\rho_7} \rangle.  $$
 This
 coincides with the \(G\)-constellation from \cite[Table 5.5, Line
 2]{Crawthesis}.
\end{example}

\subsection{Distinguished $G$-clusters}

Before specializing from $G$-constellations to $G$-clusters by
choosing $\theta\in \Theta$ so that $\mathcal{M}_{\theta}\cong\ghilb$, we prove a
pair of lemmas that are valid for more general $\theta$.

\begin{definition}
The quiver $Q_b$ associated to a representation $b=(b_i^{\rho}) \in Z$
is the subquiver of the McKay quiver with a vertex for each $\rho \in
G^*$ and those arrows $a_i^{\rho}$ for which $b_i^{\rho} \neq 0$.
When $b= b_{\theta,\mathbf{w}}$, we write $Q_{\theta,\mathbf{w}}$.
\end{definition}

In what follows we identify the support of a vector
$\mathbf{u}\in\Q^{nr}$ with the quiver containing the arrows
$a_i^{\rho}$ for which $u_i^{\rho} \neq 0$.  Conversely, given an
(undirected) path in the McKay quiver, its \emph{vector} is the
weighted sum of those $\mathbf{e}_i^{\rho}$ for which $a_i^{\rho}$
appears in the path, with the weight recording the number of times the
edge is crossed in the forward direction minus the number of times it
is crossed in the negative direction.  Recall that the lattice $M
\subset \mathbb Z^n$ is the kernel of the map $\deg \colon \mathbb Z^n
\rightarrow G^*$, and that $P_\theta$ is a polyhedron whose normal fan
is equal to the toric fan of $Y_\theta$ (see Theorem~\ref{thm:CMT1}(3)
for the construction). Note that since $Y_\theta$ is constructed by
GIT, we have $Y_{\theta} = Y_{j\theta}$ for any $j>0$.

\begin{lemma} \label{l:usupp} Fix $\theta' \in \Theta \cap \mathbb Z^r$ and 
${\bf w}\in (\Q^n_{\geq 0})^*$. Let $N=(nr)!$ and set
  $\theta=N\theta'$.  Let $F_\theta$ be the face of $P_\theta$
  minimizing ${\bf w}$.  Then there exists ${\bf u}\in \N^{nr}$ and
  ${\bf m}\in \relint(F_\theta)\cap M$ such that $C{\bf u} =
  (\theta,\mathbf{m})$ and $\supp({\bf u}) = Q_{\theta,{\bf w}}$.  In
  particular, if $\mathbf{m}$ is a vertex of $P_{\theta}$, then there
  exists $\mathbf{u} \in \mathbb N^{nr}$ such that
  $C\mathbf{u}=(\theta,\mathbf{m})$ and the quiver $\supp({\bf u}) =
  Q_{\theta,{\bf w}}$ contains no directed cycles.
\end{lemma}

\begin{proof}
  Let $F$ be the smallest face of $P$
  containing the preimage under $\pi_n$ of the face $F_{\theta'}$ of
  $P_{\theta'}$ minimizing ${\bf w}$.  This face corresponds to the
  distinguished representation $b_{\theta',\mathbf{w}}=(b_i^{\rho})$
  by definition, so $F$ is the positive rational span of those
  $C_i^{\rho}$ with $b_i^{\rho} \neq 0$ by
  Theorem~\ref{thm:main}. The fact that the matrix $B$
  is unimodular implies that there are integral points on all faces of
  $\{ {\bf u} \geq {\bf 0} \,:\, B{\bf u} = \theta' \}$.  Since
  $P_{\theta'} =\conv(\{\pi_n(C\mathbf{u}) : \mathbf{u} \in \mathbb
  Q_{\geq 0}^{nr}, B\mathbf{u}=\theta'\} )$, there is a lattice point
  $\mathbf{m}' \in \pi_n(\relint(F \cap \pi^{-1}(\theta'))) \cap M$.
  Since $(\theta',\mathbf{m}')$ lies in the relative interior of $F$
  there is $\mathbf{u}' \in \mathbb Q_{\geq 0}^{nr}$ with
  $C\mathbf{u}' =(\theta',\mathbf{m}')$ such that the support of
  $\mathbf{u}'$ is contained in $Q_{\theta,\mathbf{w}}$.  Again, since
  $B$ is unimodular we may take $\mathbf{u}' \in \mathbb N^{nr}$.
  Furthermore, for every arrow $a_i^{\rho}$ in $Q_{\theta,\mathbf{w}}$
  there is such a $\mathbf{u}'$ with $u_i^{\rho} >0$.  Adding these
  together gives a vector $\mathbf{u}'' \in \mathbb N^{nr}$ with
  support exactly $Q_{\theta,\mathbf{w}}$, and
  $C\mathbf{u}''=(f\theta',f\mathbf{m}')$, where $f$ is the number of
  arrows in $Q_{\theta,\mathbf{w}}$.  Since $f$ divides $(nr)!$ we set
  $\mathbf{u}:=(N/f)\mathbf{u}''$, which satisfies $\mathbf{u} \in
  \mathbb N^{nr}$, $\supp(\mathbf{u}) = Q_{\theta,\mathbf{w}}$, and
  $C\mathbf{u}=(\theta,N\mathbf{m}')$.  The result now follows by
  setting $\mathbf{m}=N\mathbf{m}'$, and observing that since
  $\mathbf{m}'$ lies in the face of $P_{\theta'}$ minimizing
  $\mathbf{w}$, we must have $\mathbf{m}$ in the face of $P_{\theta}$
  minimizing $\mathbf{w}$.
  
  To prove the final statement, suppose that ${\bf m}\in P_\theta$ is
  a vertex and that the quiver $Q_{\theta,{\bf w}}$ contains a
  directed cycle consisting of \(\alpha_i\) arrows labeled \(i\).  By
  adding together the collection of all equations $w_i + v_\rho =
  v_{\rho\rho_i}$ arising via (\ref{eqn:birho}) from each arrow in the
  cycle, we obtain \(\sum \alpha_i w_i = 0\).  This means that ${\bf
    w}$ is constrained to lie in a hyperplane, but this is absurd
  since ${\bf w}$ may be any vector in the relative interior of the
  top-dimensional cone $\mathcal{N}_{P_\theta}({\bf m})$ which consists of
  all ${\bf w}'\in (\Q^n_{\geq 0})^*$ such that the linear functional
  ${\bf w}'$ is minimized over $P_\theta$ at ${\bf m}$.
\end{proof}

\begin{lemma} \label{lemma:pathdecomp} Fix $\theta \in \Theta$ with
  $\theta_{\rho_0} \leq 0$ and $\theta_{\rho} \geq 0$ for $\rho \neq
  \rho_0$.  Let ${\bf u}\in \N^{nr}$ be such that $B{\bf u} = \theta$.
  Then ${\bf u}$ decomposes as ${\bf u} = {\bf u}_0 + \sum_k {\bf
    u}_k$ where ${\bf u}_0$ is the vector of a union of cycles, and, for
  each $k$, we have ${\bf u}_k\in \N^{nr}$, ${\bf u}_k\leq {\bf u}$
  and $B{\bf u}_k = {\bf e}_{\rho} - {\bf e}_{\rho_0}$ for some
  $\rho\in G^*$ depending on $k$.  Each ${\bf u}_k$ is the vector of a
  directed path in the McKay quiver from $\rho_0$ to $\rho$.
\end{lemma}

\begin{proof}
  The proof is by induction on $\sum_{\rho \in G^*} |\theta_{\rho}|$.
  When this sum is zero, $\theta=0$, so $\mathbf{u} \in \ker_{\mathbb
    Z}(B) \cap \mathbb N^{nr}$, and thus by Exercise 38 of
  Bollob{\'a}s~\cite[II.3]{Bollobas} there is a collection of directed
  cycles $\gamma_k$ in the McKay quiver with vectors $\mathbf{u}_k \in
  \mathbb N^{nr}$ with $B{\mathbf{u}_k}=0$, and $\sum_k \mathbf{u}_k =
  \mathbf{u}$ as required.  We may then assume that the lemma is true
  for smaller $\sum |\theta_{\rho}|$, and that this sum is positive,
  so $\theta_{\rho_0}<0$.

  Let $\mathcal A_{\bf u}$ be the collection of arrows in the McKay
  quiver consisting of $u_i^{\rho}$ copies of the arrow $a_i^{\rho}$
  for each pair $(i,\rho)$ with $u_i^{\rho} > 0$.  Since
  $\theta_{\rho_0}<0$, there exists some $i$, $1 \leq i \leq n$, such
  that $u_i^{\rho_0\rho_i^{-1}}=u_i^{\rho_i^{-1}} > 0$ or
  equivalently, there exists some $i$ such that $a_i^{\rho_i^{-1}} \in
  \mathcal A_{\bf u}$.  This is the first arrow in a path in the McKay
  quiver from $\rho_0$ consisting of arrows from $\mathcal{A}_{\bf
    u}$.  Continue this path until you reach a vertex $\rho$ that has
  no arrows in $\mathcal A_{\mathbf{u}}$ with tail at $\rho$.  This
  means that all the entries in the row of $B$ indexed by $\rho$ and
  lying in columns indexed by arrows in $\mathcal A_{\bf u}$ are
  $+1$s. Since $B{\bf u} = \theta$, we have $\theta_{\rho} > 0$.  This
  constructs a path $\gamma_{\rho}$ from $\rho_0$ to $\rho$ with ${\bf
    u}_1:=\mathbf{v}(\gamma_{\rho})$ using arrows in $\mathcal A_{\bf
    u}$. By construction, ${\bf u}_1 \leq {\bf u}$ and $B{\bf u}_1 =
  {\bf e}_{\rho} - {\bf e}_{\rho_0}$.  Let
  $\mathbf{u}'=\mathbf{u}-\mathbf{u}_1$.  Then
  $B\mathbf{u}'=\theta-\mathbf{e}_{\rho}+\mathbf{e}_{\rho_0}$, which
  has smaller coordinate sum, so by the induction hypothesis
  $\mathbf{u}'$, and thus $\mathbf{u}$, has a decomposition of the
  desired form.
\end{proof} 

For the rest of this paper we restrict to $G$-clusters. To do this, fix a
 parameter $\theta\in \Theta$ satisfying \begin{equation}
   \label{eqn:ghilbtheta} \left\{\begin{array}{l} \theta_{\rho_0} < 0
       \mbox{ and }\theta_\rho > 0 \mbox{ for all }\rho\neq \rho_0; \text{
         and } \\ \theta/(nr)! \in \Z^{r} \text{ has }\oplus_{j \geq 0}
       \K[V]_{j\theta/(nr)!} \text{ generated in degree
         one}.\end{array}\right.  \end{equation}
 The result of Ito--Nakajima~\cite[\S2]{ItoNakajima} implies that
 $\theta$-stable $G$-constellations are precisely $G$-equivariant
 $S$-modules of the form \(S/J\), where \(J\subseteq S\) is a
 $G$-invariant ideal and \(S/J\) is isomorphic to $\K G$ as a \(\K
 G\)-module.  Thus, $\mathcal{M}_\theta \cong \ghilb$, and hence
 $Y_\theta \cong \hilbg$.  The second assumption in
 (\ref{eqn:ghilbtheta}) guarantees that \(\theta' := \theta/(nr)!\)
 lies in $\mathbb Z^r$, so we can apply Lemma~\ref{l:usupp}.  This
 second assumption is required only for the proofs and is not relevant
 when computing examples.

\begin{definition}
 \label{defn:IM}
Let $I_M$ be the lattice ideal $\langle x^{\mathbf{u}}
-x^{\mathbf{u}'} : \mathbf{u}, \mathbf{u}' \in \mathbb N^n,
\mathbf{u}-\mathbf{u}' \in M \rangle$. The scheme $Z(I_M)\subseteq
\mathbb{A}^n_{\K}$ is the $G$-orbit of the point $(1,\dots,1)\in
\mathbb{A}^n_{\K}$.
\end{definition}

Lattice ideals are generalizations of toric ideals, and have many
applications.  See, for example, the book of
Miller-Sturmfels~\cite{MillerSturmfels}.

\begin{proposition}
  \label{prop:vertexideal} Let $J \subseteq S$ be an ideal
  defining a point \([J]\in \ghilb\).  Then \([J]\) defines the
  distinguished point \([b_{\theta,\mathbf{w}}]\in \hilbg\) if and
  only if $J=\inn_{\mathbf{w}}(I_M)$.
\end{proposition} 

\begin{proof} 
  The ideal \(J\subseteq S\) defines a point \([J]\in \hilbg\) if and
  only if \([J] = [b_{\theta,\mathbf{w}}]\) for some \({\bf w}\in
  (\Q^n_{\geq 0})^*\), in which case \(J\) satisfies $S/J\cong
  A/M_{\theta,\mathbf{w}}$ as an $S$-module.
  
  We claim that \(J\) is the kernel of the \(S\)-module homomorphism
  $\phi \colon S \rightarrow A/M_{\theta,\mathbf{w}}$ defined by
  setting $\phi(1)=\mathbf{e}_{\rho_0}$.  Indeed, by
  Lemma~\ref{l:usupp} there is a vector $\mathbf{u} \in \mathbb
  N^{nr}$ with support $Q_{\theta,\mathbf{w}}$ and
  $B\mathbf{u}=\theta$.  By Lemma~\ref{lemma:pathdecomp} we may
  decompose $\mathbf{u}$ as a sum of vectors $\mathbf{u}_k \leq
  \mathbf{u}$, with ${\bf u}_k$ the vector of a path from $\rho_0$ to
  some $\rho\in G^*$, and a vector ${\bf u}_0 \leq {\bf u}$ with ${\bf
    u}_0$ the vector of a union of cycles.  Since the support of
  $\mathbf{u}$ is $Q_{\theta,\mathbf{w}}$, these paths and cycles lie
  in $Q_{\theta,\mathbf{w}}$.  By assumption (\ref{eqn:ghilbtheta}) on
  $\theta\in \Theta$, there is at least one such vector ${\bf u}_k$
  for each $\rho \in G^*$.  This implies that $\phi$ is surjective, as
  a path from $\rho_0$ to $\rho$ in $Q_{\theta,\mathbf{w}}$ yields a
  binomial of the form $x^{\mathbf{v}}\mathbf{e}_{\rho_0}-\lambda
  \mathbf{e}_{\mathbf{\rho}}$ in $M_{\theta,\mathbf{w}}$ with $\lambda
  \neq 0$.  This gives $S/\ker(\phi)\cong A/M_{\theta,\mathbf{w}}$, so
  \(J = \ker(\phi)\) as required.
  
  It remains to show that $\ker(\phi)=\inn_{\mathbf{w}}(I_M)$.  The
  \(G^*\)-graded Hilbert function of \(S/I_M\) is one in every degree
  by Definition~\ref{defn:IM}, and thus the same is true of
  \(S/\inn_{\bf w}(I_M)\) for \({\bf w}\in (\Q^n_{\geq 0})^*\).  Since
  $S/\ker(\phi)\cong A/M_{\theta,\mathbf{w}}$, it follows that
  \(\ker(\phi)\) and $\inn_{\mathbf{w}}(I_M)$ have the same Hilbert
  function. If $x^{\mathbf{u}}$ is a minimal generator of
  $\inn_{\mathbf{w}}(I_M)$, then since $S/\inn_{\mathbf{w}}(I_M)$ has
  $G^*$-graded Hilbert function one in every degree, there is
  $x^{\mathbf{u}'} \not \in \inn_{{\bf w}}(I_M)$ with
  $\mathbf{u}-\mathbf{u}' \in M$ and $\mathbf{w} \cdot \mathbf{u} >
  \mathbf{w} \cdot \mathbf{u}'$.  Now
  $(x^{\mathbf{u}}-x^{\mathbf{u}'})\mathbf{e}_{\rho_0} \in M_G$ since
  this binomial is homogeneous under $G^{\ast}$-grading, so
  $x^{\mathbf{u}}\mathbf{e}_{\rho_0} \in
  \inn_{(\mathbf{v},\mathbf{w})}(M_G)$ for any $\mathbf{v} \in
  (\mathbb Q^r)^*$, and thus $x^{\mathbf{u}}\mathbf{e}_{\rho_0} \in
  M_{\theta,\mathbf{w}}$.  The proof is identical for a minimal
  generator $x^{\mathbf{u}}-x^{\mathbf{u}'}$ of
  $\inn_{\mathbf{w}}(I_M)$, so $\inn_{\mathbf{w}}(I_M)\subseteq
  \ker(\phi)$.  Since $\ker(\phi)$ and $\inn_{\mathbf{w}}(I_M)$ have
  the same Hilbert function, they must in fact be equal.  This
  completes the proof.  \end{proof}

If $\gamma$ is a path in the McKay quiver, its \emph{type} is the vector $\mathbf{u} =(u_i) \in \mathbb Z^n$ with $u_i$ being the number of arrows labelled $i$ in $\gamma$.

\begin{corollary} \label{c:quiverJ} Let ${\bf w}\in (\Q^{n}_{\geq
    0})^*$ lie in the relative interior of a top-dimensional cone in
  the fan of $Y_\theta$. Then a directed path from $\rho_0$ to $\rho$
  of type ${\bf u}\in \N^n$ is supported on the arrows in
  $Q_{\theta,{\bf w}}$ if and only if $x^{\bf u}\not\in \inn_{{\bf
      w}}(I_M)$ and $\deg(x^{\bf u}) = \rho^{-1}$.
 \end{corollary}

 We now exhibit an example to illustrate that \(\ghilb\)
 need not be irreducible, thereby proving that
 \(\mathcal{M}_\theta\neq Y_\theta\) in general.  It is known by
 Bridgeland--King--Reid~\cite{BKR} and Ishii~\cite{Ishii1} that
 \(\ghilb\) is irreducible for finite subgroups of \(\SL(n,\K)\) with
 \(n\leq 3\) and \(\GL(2,\K)\) respectively.  Thus, the simplest
 possible reducible example is determined by a finite subgroup of
 \(\GL(3,\K)\).  

 \begin{example}
 \label{ex:reducible}
 Consider the group \(G:= \Z/14\Z\) embedded in \(\GL(3,\K)\) with
 generator the diagonal matrix \(g =
 \diag(\omega^{1},\omega^9,\omega^{11})\), where \(\omega\) is a
 primitive fourteenth root of unity.  We claim that for \(\theta\in
 \Theta\) as in (\ref{eqn:ghilbtheta}), the moduli space
 \(\mathcal{M}_{\theta} \cong \ghilb\) is reducible, hence \(\ghilb \neq
 \hilbg\).  To see this, consider the torus-fixed point
 \([J]\in\ghilb\) defined by the monomial ideal
 \[
 J = \langle x_2^2x_3,
 x_1x_3^2, x_1x_2^2, x_1^2x_2,x_2x_3^2, x_1^2x_3, x_2^4, x_3^4, x_1^4 \rangle
 \]
 in \(\K[x_1,x_2,x_3]\).  We check that $J$ is a $G$-cluster by
 checking its $G^*$-graded Hilbert function.  In this example this can
 be done by hand or using Macaulay 2 \cite{M2}.
 
 To show that \([J]\) does not lie in \(\hilbg\) we establish that $J$
 is not an initial ideal of \(I_M=\langle x_1^{14}-1, x_2-x_1^{9},
 x_3-x_1^{11}\rangle\).  Suppose otherwise, so there is a weight
 vector \({\textbf w}\in (\Q^3_{\geq 0})^*\) with \({\textbf w} \cdot
 {\textbf u} > {\textbf w} \cdot {\textbf u'}\) whenever
 \({x}^{\textbf u}-{x}^{\textbf u'} \in I_M\) for \({x}^{\textbf u}
 \in J\) and \({x}^{\textbf u'} \notin J\).  Consider the binomials
 \({\underline {x_1^{2}x_3}}-x_2^3\), \({\underline {x_2x_3^{2}}}-x_1^3\)
 and \({\underline {x_1x_2^2}}-x_3^3\) in \(I_M\) where the underlined
 monomials are minimal generators of $J$ and the trailing monomial in
 each binomial is the unique standard monomial of $J$ in its degree.
 If $J$ was $\textup{in}_{\bf w}(I_M)$ for a weight vector
 \({\textbf w} = (w_1,w_2,w_3)\) then these binomials would 
imply that
\[
2w_{1}+w_{3} > 3w_{2}, \quad w_{2} + 2w_{3} > 3w_{1}, \quad w_{1} +
2w_{2} > 3w_{3}.
 \]
 Adding these three inequalities leads to the new inequality
 \(3w_1+3w_2+3w_3 > 3w_1+3w_2+3w_3\), which is absurd.
\end{example}

  \begin{remark}
 \label{remark:counterexample}
 The monomial ideal \(J\) constructed in
 Example~\ref{ex:reducible} defines a point \([J]\in \ghilb\) that
 lies off the coherent component \(\hilbg\).  This provides a
 counterexample to the statements of Nakamura~\cite[Corollary~2.4,
 Theorem~2.11]{Nakamura} that every monomial \(G\)-cluster \(J\subset
 S\) defines a point \([J]\in \hilbg\). For each monomial $G$-cluster $J$, Nakamura~\cite[\S 1]{Nakamura}
 defined a cone $\sigma(\Gamma)$ associated to the set of standard
 monomials $\Gamma$ (which he calls a \emph{$G$-graph}).  If
 $J=\inn_{\mathbf{w}}(I_M)$, this is the cone in the Gr\"obner fan of
 $I_M$ corresponding to $J$, while this set is empty otherwise.  See
 \cite[Chapter 2]{GBCP} or \cite[Chapter 2]{IndiaNotes} for details on
 the Gr\"obner fan.
\end{remark}

\section{Local equations on Nakamura's \protect$G\protect$-Hilbert scheme}
\label{sec:nonnormal}

In this section we give a new description of local coordinate charts
on $\hilbg$.  We use this description to exhibit a finite subgroup $G
\subset \GL(6,\K)$ for which $\hilbg$ is nonnormal.

 \subsection{Local equations on $\hilbg$}
 In \cite{CMT1}, we showed that when $\theta$ satisfies the conditions
 of (\ref{eqn:ghilbtheta}), local charts on the coherent
 component $Y_\theta \cong \hilbg$ are given by $\Spec(\K[A_{\sigma}])$,
 where $\sigma=\mathcal N_{P_{\theta}}(\mathbf{m}):= \{{\bf w}'\in
 (\Q^n)^* : {\bf w}' \text{ is minimized over }P_\theta \text{ at
 }{\bf m}\}$ for a vertex
 $\mathbf{m}$ of $P_{\theta}$, and
 \[
 A_\sigma = \mathbb N\langle \mathbf{p}-\mathbf{m} : \mathbf{p} \in P_\theta\cap
M \rangle.
 \]
 We now provide an alternative description of these local charts.
 Theorem~\ref{thm:coverhilbg} below corrects and refines the result of
 Nakamura (compare Remark~\ref{remark:counterexample}).

\begin{definition}
  Let $J=\inn_{\mathbf{w}}(I_M)$ be a monomial initial ideal.  We
  associate to $J$ the semigroup $A_{J}$ generated by
  $\{\mathbf{u}-\mathbf{u}' \in M : \mathbf{u}, \mathbf{u}' \in
  \mathbb N^{n}, x^{\mathbf{u}} \in J, x^{\mathbf{u}'} \not \in J \}.$
\end{definition}

Note that the initial ideal $\inn_{\mathbf{w}}(I_M)$ is a monomial
ideal when ${\bf w}$ is generic.  For a vector $\mathbf{u} \in \mathbb
Z^n$ we write $\pos(\mathbf{u})$ for the vector with $i$th component
$u_i$ if $u_i>0$ and $0$ otherwise.  Similarly, we write
$\negg(\mathbf{u})$ for the vector with $i$th component $-u_i$ if
$u_i<0$ and $0$ otherwise, so $\mathbf{u}=\pos(\mathbf{u}) -
\negg(\mathbf{u})$.  Note that for ${\bf u}\in \Z^n$, we have $\deg({\bf u})\in G^*$ and $\deg(-{\bf u}) = \deg({\bf u})^{-1}\in G^*$.

\begin{theorem} \label{thm:coverhilbg} For $\theta \in \Theta$ as in
  (\ref{eqn:ghilbtheta}), let $\mathbf{m}\in P_\theta\cap M$ be a vertex and
  choose $\mathbf{w}$ in the relative interior of $\sigma =\mathcal
  N_{P_{\theta}}(\mathbf{m})$.  Then $A_{\sigma} = A_{J}$, where
  $J=\inn_{\mathbf{w}}(I_M)$.  Thus \(\hilbg\) is covered
  by affine charts \(\Spec \K[A_J]\) defined by the monomial ideals
  \(J=\inn_{\mathbf{w}}(I_M)\) as ${\bf w}$ varies in $(\mathbb
  Q^n_{\geq 0})^{\ast}$.
 \end{theorem}

 \begin{proof} 
   We first show that $A_{J} \subseteq A_{\sigma}$.  Let
   $\mathbf{u}-\mathbf{u'}\in A_{J}$, with $x^\mathbf{u} \in J$ and
   $x^\mathbf{u'} \not \in J$.  It is enough to establish that the
   lattice point ${\bf p}:= \mathbf{u}-\mathbf{u'}+\mathbf{m}\in M$
   lies in $P_\theta$.  To show this we construct vectors
   $\textbf{u}_\mathbf{u},\textbf{u}_\mathbf{u'},\textbf{u}_\mathbf{m}\in
   \N^{nr}$ such that $\pi_n(C(\textbf{u}_\mathbf{u} -
   \textbf{u}_\mathbf{u'} + \textbf{u}_\mathbf{m}))$ lies in $P_\theta$
   by construction and is equal to $\mathbf{p}$.

   By Lemma~\ref{l:usupp} there is a vector $\textbf{u}_\mathbf{m}\in \N^{nr}$
   satisfying $C\textbf{u}_\mathbf{m} = (\theta,\mathbf{m})\in \N C$
   with $\supp(\mathbf{u}_\mathbf{m})=Q_{\theta,\mathbf{w}}$. Next, by
   Corollary~\ref{c:quiverJ},
   since $x^\mathbf{u'}\not\in J$, any directed path
   $\gamma_\mathbf{u'}$ in the McKay quiver from $\rho_0$ to
   $\deg(-\mathbf{u'}) \in G^*$ of type $\mathbf{u'}\in\mathbb N^{n}$ is
   supported on arrows in the subquiver $Q_{\theta,\mathbf{w}}$. This
   path determines a vector $\textbf{u}_\mathbf{u'}:=
   \mathbf{v}(\gamma_\mathbf{u'})\in \N^{nr}$ satisfying
   $C\mathbf{u}_\mathbf{u'} = (\textbf{e}_{\deg(-\mathbf{u'})} -
   \textbf{e}_{\rho_0},\mathbf{u'})\in \N C$ and
   $\mathbf{u}_\mathbf{m} - \mathbf{u}_\mathbf{u'}\in \N^{nr}$.
   Finally, pick any path $\gamma_\mathbf{u}$ of type
   $\mathbf{u}\in\mathbb N^{n}$ beginning at $\rho_0$.  The head of
   this path is $\deg(-\mathbf{u})$, so $\mathbf{u}_{{\bf
       u}}:=\mathbf{\mathbf{v}}(\gamma_\mathbf{u})$ satisfies
   $C\mathbf{u}_\mathbf{u}=(\mathbf{e}_{\deg(-{\bf
       u})}-\mathbf{e}_{\rho_0},{\bf u})$.  By adding
   $\textbf{u}_\mathbf{u}\in \N^{nr}$ to $\mathbf{u}_\mathbf{m} -
   \mathbf{u}_\mathbf{u'}\in \N^{nr}$, we obtain a vector in $\N^{nr}$
   satisfying $C(\mathbf{u}_\mathbf{u} - \mathbf{u}_\mathbf{u'} +
   \mathbf{u}_\mathbf{m}) =
   (\theta,\mathbf{u}-\mathbf{u'}+\mathbf{m})$, since
   $\deg(\mathbf{u})=\deg(\mathbf{u'})$.  Hence $\mathbf{p} =
   \pi_n(C(\mathbf{u}_\mathbf{u} - \mathbf{u}_\mathbf{u'} +
   \mathbf{u}_\mathbf{m}))$ lies in $P_\theta$ as claimed.

   For the opposite inclusion, consider a minimal generator
   $\mathbf{p}-\mathbf{m}\in A_{\sigma}$. By Lemma~\ref{l:usupp},
   there exists $\textbf{u}_\mathbf{m}\in \N^{nr}$ such that $C
   \textbf{u}_\mathbf{m} = (\theta,{\bf m})$ and the quiver
   $\supp({\bf u}_\mathbf{m}) = Q_{\theta,{\bf w}}$ contains no directed cycles. Since ${\bf p}\in P_\theta$ and ${\bf
   p}\neq {\bf m}$, there exists $\textbf{u}_{\bf p}\in \N^{nr}$ such
   that $C \textbf{u}_{\bf p} = (\theta,{\bf p})$ and
   $\textbf{u}_\mathbf{m}\neq \textbf{u}_{\bf p}$.
   Lemma~\ref{lemma:pathdecomp} enables us to decompose
   $\textbf{u}_\mathbf{p}$ into a sum of vectors of the form
   $\textbf{u}_\mathbf{p}(\rho)\in \N^{nr}$, where each
   $\textbf{u}_\mathbf{p}(\rho)$ satisfies
   $B\textbf{u}_\mathbf{p}(\rho) = \textbf{e}_{\rho} -
   \textbf{e}_{\rho_0}$ and $\textbf{u}_\mathbf{p}(\rho)\leq
   \textbf{u}_\mathbf{p}$ by construction, and also a vector ${\bf
     u}_{\bf p}(0)$, where ${\bf u}_{\bf p}(0) \leq {\bf u}_{\bf
     p}$, and ${\bf u}_{\bf p}(0)$ is the vector of a union of cycles.  The same is true for
   $\textbf{u}_{\bf m}$.  For each $\rho \neq \rho_0$, there are
   $\theta_{\rho}$ vectors of the form $\mathbf{u}_\mathbf{m}(\rho)$
   and $\theta_{\rho}$ of the form $\mathbf{u}_{\bf p}(\rho)$.  Note
   that ${\bf u}_{\bf m}(0) = 0$ since $Q_{\theta,{\bf w}}$
   contains no cycles.

   There are now two cases.  Either there exists $\rho \neq \rho_0$
   and vectors $\mathbf{u}_{\bf p}(\rho), \mathbf{u}_\mathbf{m}(\rho)$
   satisfying $\pi_n(C\mathbf{u}_{\bf p}(\rho))\neq
   \pi_n(C\mathbf{u}_\mathbf{m}(\rho))$, or else $\mathbf{u}_{\bf
     p}-\mathbf{u}_\mathbf{m}= {\bf u}_{\bf p}(\rho_0)$.  In the
   latter case, ${\bf p}-{\bf m} = \pi_n(C({\bf u}_{\bf p}(\rho_0)))$
   lies in $\N^n$ and satisfies $\deg({\bf p}-{\bf m}) = \rho_0$.
   Since the only standard monomial of $J$ of degree $\rho_0$ is $1$,
   we have $x^{{\bf p}-{\bf m}} \in J$, so ${\bf p}-{\bf m} \in A_J$
   as required. In the former case, suppose that
   $\pi_n(C\mathbf{u}_{\bf p}(\rho)) \neq
   \pi_n(C\mathbf{u}_\mathbf{m}(\rho))$.  Let
   $\mathbf{u}'_\mathbf{p}=\mathbf{u}_\mathbf{p}-\mathbf{u}_\mathbf{p}(\rho)+\mathbf{u}_\mathbf{m}(\rho)$,
   and
   $\mathbf{u}'_\mathbf{m}=\mathbf{u}_\mathbf{m}-\mathbf{u}_\mathbf{m}(\rho)+\mathbf{u}_\mathbf{p}(\rho)$.
   Note that $\mathbf{u}'_\mathbf{p}, \mathbf{u}'_\mathbf{m} \in
   \mathbb N^{nr}$ and
   $B\mathbf{u}'_\mathbf{p}=B\mathbf{u}'_\mathbf{m}=\theta$, so
   $\pi_n(C\mathbf{u}'_\mathbf{p})-\mathbf{m}$ and
   $\pi_n(C\mathbf{u}'_\mathbf{m})-\mathbf{m}$ both lie in
   $A_{\sigma}$.  In addition we have
   $(\pi_n(C\mathbf{u}'_\mathbf{m})-\mathbf{m})+(\pi_n(C\mathbf{u}'_\mathbf{p})-\mathbf{m})=\pi_n(C(\mathbf{u}'_\mathbf{m}-\mathbf{u}_\mathbf{m}
   +\mathbf{u}'_\mathbf{p}-\mathbf{u}_\mathbf{m}))=\pi_n(C(\mathbf{u}_\mathbf{p}-\mathbf{u}_\mathbf{m}))=\mathbf{p}-\mathbf{m}$.
   Since $\mathbf{p}-\mathbf{m}$ is a minimal generator of
   $A_{\sigma}$, and
   $\pi_n(C(\mathbf{u}_\mathbf{p}(\rho)-\mathbf{u}_\mathbf{m}(\rho)))\neq
   0$, the second of these terms $\pi_n(C{\bf u}'_{\bf p}) - {\bf m}$
   must be zero.  This gives $\mathbf{p}-\mathbf{m} =
   \pi_n(C(\mathbf{u}_\mathbf{p}(\rho)-\mathbf{u}_\mathbf{m}(\rho)))$,
   where $\mathbf{u}_\mathbf{p}(\rho)$ and
   $\mathbf{u}_\mathbf{m}(\rho)$ are vectors of paths from $\rho_0$ to
   $\rho$.  Since \(\mathbf{u}_{\bf m}(\rho)\leq {\bf u}_{\bf m}\),
   the support of the path defined by $\mathbf{u}_\mathbf{m}(\rho)$
   lies in \(Q_{\theta, {\bf w}}\), so
   $x^{\pi_n(C(\textbf{u}_\mathbf{m}(\rho))}\not\in J$ by
   Corollary~\ref{c:quiverJ}.  Also, the negative part of ${\bf
     p}-{\bf m}$ satisfies $\negg({\bf p}-{\bf m})\leq
   \pi_n(C\negg(\textbf{u}_{\bf p}(\rho) -
   \textbf{u}_\mathbf{m}(\rho))) \leq
   \pi_n(C\textbf{u}_\mathbf{m}(\rho))$.  This gives
   $x^{\negg({\bf p}-{\bf m})} \not \in J$.  Since ${\bf p}-{\bf m}
   \neq 0$, and $J$ has only one standard monomial of each degree, we
   must have $x^{\pos({\bf p}-{\bf m})} \in J$, so ${\bf p}-{\bf m}
   \in A_J$ as required.  \end{proof}

We next give a smaller generating set for $A_J$.

\begin{lemma}
 \label{l:naklemma}
  The semigroup $A_J$ is generated by elements of the form $
  \{\mathbf{u}-\mathbf{u}' \in M : \mathbf{u},\mathbf{u}'\in \N^n,
  \deg(\mathbf{u}) = \deg(\mathbf{u}'), x^{\mathbf{u}} \text{ is a
    minimal generator of }J, x^{\mathbf{u}'}\not\in J\}$.
\end{lemma}

\begin{proof}
  The semigroup $A'$ generated by the given elements is a subsemigroup
  of $A_J$, so we need only show that if $\mathbf{u}-\mathbf{v} \in
  A_J$, with $x^\mathbf{u} \in J$, $x^\mathbf{v} \not \in J$, and
  $\deg(x^\mathbf{u})=\deg(x^\mathbf{v})$, then
  $\mathbf{u}-\mathbf{v}$ is in the semigroup generated by $A'$. Note
  that $\mathcal G=\{x^\mathbf{u}-x^\mathbf{v} : \mathbf{u}-
  \mathbf{v} \in A'\}$ is a Gr\"obner basis for $I_M$.  Since
  $x^{\mathbf{v}}$ is the unique standard monomial of $J$ of its
  degree, this means that $x^\mathbf{u}$ reduces modulo $\mathcal G$
  to $x^\mathbf{v}$.  So we can write
  $x^{\mathbf{u}}-x^{\mathbf{v}}=\sum_{j=1}^s
  x^{\mathbf{w}_{i_j}}(x^{\mathbf{u}_{i_j}}-x^{\mathbf{v}_{i_j}})$,
  where $x^{\mathbf{u}_{i_j}}-x^{\mathbf{v}_{i_j}} \in \mathcal G$,
  $x^{\mathbf{w}_{i_j}}x^{\mathbf{u}_{i_j}}=x^{\mathbf{w}_{i_{j-1}}}x^{\mathbf{v}_{i_{j-1}}}$
  for $2 \leq j \leq s$, and
  $x^{\mathbf{w}_{i_1}}x^{\mathbf{u}_{i_1}}=x^{\mathbf{u}}$,
  $x^{\mathbf{w}_{i_s}}x^{\mathbf{v}_{i_s}}=x^{\mathbf{v}}$.  This
  means that $\mathbf{u}-\mathbf{v}=\sum_{j=1}^s
  \mathbf{u}_{i_j}-\mathbf{v}_{i_j}$, so $\mathbf{u}-\mathbf{v}$ lies
  in the semigroup generated by $A'$.
\end{proof}
 
 \begin{remark}
   Lemma~\ref{l:naklemma} is the content of
   Nakamura~\cite[Lemma~1.8]{Nakamura}.  We provide a self-contained
   proof to illustrate the Gr\"obner argument.  Note that
   Example~\ref{ex:reducible} is a counterexample to the sentence
   following \cite[Lemma~1.8]{Nakamura}.
 \end{remark}

 List the elements from the generating set of $A_J$ presented in
Lemma~\ref{l:naklemma} as $\{\mathbf{u}_1-\mathbf{u}'_1, \dots,
\mathbf{u}_s-\mathbf{u}'_s\}$. Let $I_U$ denote the kernel of the
$\K$-algebra homomorphism $\K[y_1,\dots,y_s]\to \K[A_J]$ sending $y_i$
to $x^{\mathbf{u}_i}/x^{\mathbf{u}'_i}$. This ideal defines the local chart
$U:= \Spec \K[A_J]$ in $\hilbg$.

 \begin{corollary} The universal family above the chart $\Spec\K[A_J]$ is given by 
 $$
 F:=\langle x^{\mathbf{u}_i} - y_i x^{\mathbf{u}'_i} : 1\leq i\leq s \rangle + I_U
 $$
in the ring $\K[x_1,\dots,x_n][y_1,\dots,y_s].$
 \end{corollary}

\begin{proof}
  Write $R=\K[x_1,\dots,x_n][y_1,\dots,y_s]$.  Let
  $\mathcal{Z}_U:=\Spec(R/F)$.  We must show that the map
  $\mathcal{Z}_U \rightarrow U$ is a flat family of $\K$-schemes with
  $\mathcal{Z}_U$ being $G$-invariant, and that $H^0(\mathcal
  O_{\mathcal{Z}_u}) \cong \K[G]$ for all $u \in U$ where
  $\mathcal{Z}_u$ is the geometric fiber over $u$.

  We first exhibit a Gr\"obner basis for $F$.  Let $\mathbf{w} \in
  (\mathbb Q^n_{\geq 0})^{\ast}$ be a weight vector for which
  $J=\inn_{\mathbf{w}}(I_M)$.  We extend $\mathbf{w}$ to
  $\tilde{\mathbf{w}} \in (\mathbb Q^{n+s}_{\geq 0})^{\ast}$ by setting
$\tilde{\mathbf{w}}_i=\mathbf{w}_i$ for $1 \leq i \leq n$, and
$\tilde{\mathbf{w}}_i=0$ for $i>n$.  Let $\prec_{\tilde{\mathbf{w}}}$
be the term order on $R$ given by refining the order given by
$\tilde{\mathbf{w}}$ by the lexicographic order.  We claim that
$\{x^{\mathbf{u}_i} - y_i x^{\mathbf{u}'_i} : 1\leq i\leq s\} \cup
\mathcal G$ is a Gr\"obner basis for $F$ with respect to
$\prec_{\tilde{\mathbf{w}}}$, where $\mathcal G$ is a Gr\"obner basis
for $I_U$ in the lexicographic order.  Indeed, since $F$ is a binomial
ideal with coefficients $\pm 1$, Buchberger's algorithm ensures that
the reduced Gr\"obner basis for $F$ with respect to
$\prec_{\tilde{\mathbf{w}}}$ also consists of binomials of this form.
Note also that if we set $\deg(x^{\mathbf{u}})=\mathbf{u} \in \mathbb
Z^n$, and $\deg(y_i)=\mathbf{u}_i-\mathbf{u}'_i$ then $F$ is
homogeneous with respect to this $\mathbb Z^n$-grading, and that this
refines the $G^*$-grading given by setting $\deg(y_i)=0$.

  Let $y^{\alpha} x^{\beta} - y^{\gamma} x^{\delta}$ be a homogeneous
  binomial in $F$ under the $\mathbb Z^n$-grading with $\beta \neq
  \delta$.  Then since $\deg(x^{\beta}) = \deg(x^{\delta})$ either
  $x^{\beta}$ or $x^{\delta}$ lies in $J$ and hence this binomial is
  reducible by an element in the first part of our proposed Gr\"obner
  basis.  Thus if the given set is not a Gr\"obner basis, then there
  exists an element of the form $(y^{\alpha} - y^{\gamma}) x^{\beta}$
  in the true Gr\"obner basis where we may have $\gamma = {\bf 0}$.
  Then since $\deg(y^{\alpha}) = \deg(y^{\gamma})$ under the $\mathbb
  Z^n$-grading, $y^{\alpha} - y^{\gamma} \in I_U$.  But then it can be
  reduced to zero using the Gr\"obner basis $\mathcal G$.  Thus, no
  such binomial exists, so the given set is a Gr\"obner basis for $F$.

This Gr\"obner basis means that $\{x^\mathbf{u} : x^\mathbf{u} \not
\in J\} \subset R$ is a basis for $R/F$ as a
$\K[y_1,\dots,y_s]/I_U$-module, so $R/F$ is a free
$\K[y_1,\dots,y_s]/I_U=\K[A_J]$-module. This implies that the map
$\mathcal{Z}_U \rightarrow U$ is flat.  Since $F$ is homogeneous with
respect to the $G^*$-grading, where $\deg(y_i)=0$ for $1\leq i\leq s$,
the scheme $\mathcal{Z}_U$ is $G$-invariant.  Since $\Spec(\K[A_J])$
is reduced, the fiber at a point $u = (u_1,\dots,u_s) \in U$ is
obtained by specializing the values of the $y_i$, and thus $\mathcal
O_{\mathcal{Z}_u} = \K[x_1,\dots,x_n]/F_u$, where $F_u$ is the result
of specialization.  The Gr\"obner result implies that $J=\inn_{\mathbf{w}}(F_u)$,
so $\K[x_1,\dots,x_n]/F_u$ has the same $G^*$-graded Hilbert function
as $J$.  By Corollary~\ref{c:gconisgraded} we conclude that
$\K[x_1,\dots,x_n]/F_u$ is a $G$-constellation, so $H^0(\mathcal
O_{\mathcal{Z}_u}) \cong \K[G]$.
\end{proof}

 \subsection{An example of a nonnormal $G$-Hilbert scheme} Recall that
 a subsemigroup $\mathbb N E$ of $\mathbb Z^n$ is {\em normal} (or
 {\em saturated}) if $ \mathbb N E= \mathbb Z E \cap \mathbb Q_{\geq
 0} E$, and that a semigroup algebra is normal as a $\K$-algebra if
 and only if the corresponding semigroup is normal.
 Theorem~\ref{thm:coverhilbg} implies that if $A_{J}$ is not a normal
 semigroup for some $J=\inn_{\mathbf{w}}(I_M)$ then the toric variety
 $\hilbg$ is not normal.

\begin{algorithm} \label{a:nonnormal}
{To check whether $\hilbg$ is normal for a given $G \subseteq \GL(n, \K)$.}

\medskip

\noindent{\bf Input: } A generating set $\mathcal L$ for the lattice
$M=\ker_\Z(\deg) \subset \mathbb Z^n$.

\begin{enumerate}
\item Compute the {\em lattice ideal} $I_{M} := \langle x^{\bf u} -
x^{\bf v} \, : \, {\bf u} - {\bf v} \in M, {\bf u}, {\bf v} \in
\mathbb N^n \rangle$. To do this, we use the result of
Ho\c{s}ten-Sturmfels~\cite{HostenSturmfels} that
  $$I_M = \bigl( \langle x^{\bf u} - x^{\bf v} \, : \, {\bf u} - {\bf
    v} \in \mathcal L, {\bf u}, {\bf v} \in \mathbb N^n \rangle \,:\,
  (\prod x_i)^\infty \bigr).$$
\item Compute all reduced Gr\"obner bases of $I_M$. This computation
  can be done using the software package Gfan \cite{gfan}.
\item For each reduced Gr\"obner basis $\mathcal G = \{x^{{\alpha}_i}
- x^{\beta_i}, \, i = 1, \ldots, t \}$, check whether the semigroup
$\mathbb N \{ \alpha_i - \beta_i, \, i = 1, \ldots, t \}$ is normal.
This can be done using the software package Normaliz \cite{Normaliz}.
 If all semigroups checked above are normal, then $\hilbg$ is
  normal. 
\end{enumerate}
\end{algorithm}

 \begin{example}
 \label{ex:nonnormal}
 Let $G\subset \GL(6,\K)$ be the subgroup generated by the diagonal
 matrices $\diag(\omega, \omega, \omega,
 \omega,\omega, \omega)$, $\diag(1, \omega, 1,
 \omega^3, \omega^4, \omega^3)$, $\diag(\omega^3,
 \omega^2, \omega^4, \omega^2, \omega,
 \omega)$, and  $\diag(\omega,1, \omega, 1, 1, 1)$, where
 $\omega$ is a primitive fifth root of unity.  The group \(G\) is
 isomorphic to \((\mathbb Z/5 \mathbb Z)^4\).  Indeed, all four
 generators have order five, and the matrix
 \[
 \begin{pmatrix}
 1&1&1&1&1&1\\
 0&1&0&3&4&3\\
 3&2&4&2&1&1\\
 1&0&1&0&0&0\\
 \end{pmatrix}
 \]
 with entries in $\mathbb Z/5 \mathbb Z$ has rank $4$, so no
 generators are redundant.  The ideal $I_M \subset \K[a,b,c,d,e,f]$ is
 $$
 I_M=\langle f^5-1, e^2f^3-b^4d, e^4f-b^3d^2, e^5-1, df^2-be^2,
 de^3-bf^3, d^2e-b^2f, d^3-b^3ef^4, $$
 $$
 ce^3-abf^2, cde-ab^2,
 c^2ef-a^2b^2, c^3f^4-a^3b^3e, c^3d^2-a^3f^2,c^4d-a^4f,c^5-1, $$
 $$
 be^2f^2-a^4c, bde^2f-a^3c^2, bc^2d^2-a^2e^3, b^2e^4-a^3c^2f,
 b^2cd^2-aef^3, b^3de-ac^4, $$
 $$
 b^4e^3-ac^4f^2, b^4cf^3-ae^2,
 b^5-1,ae^4-b^3cd, ad-cf,ac^4e^2-b^4f^3, ae^2f^2-b^4c, $$
 $$
 abe^2-cf^3,abc^4f^2-e^3, ab^2c^4-de, ab^3ef^3-cd^2, a^2f^4-b^3c^2e,
 a^2ef^2-b^2c^2d, $$
 $$
 a^2c^3-b^3ef,a^3ef-b^2c^3,a^3e^3-bc^3df,a^3c^2f^2-d^2,
 a^3c^2e-b^2f^4, a^3b^3c^2f-e^4, $$
 $$
 a^3b^4c^2-de^2f,
 a^4f^3-bc^4e^2, a^4cf-d, a^4ce^2-b^4d^2f, a^4b^4c-e^2f^2, a^5-1
 \rangle $$
 We claim that the ideal $$
 J=\langle f^5, e^2f^3, e^4f,
 e^5, df^2, de^3, d^2e, d^3, ce^3, cde, c^2ef, c^3f^4, c^3d^2, c^4d,
 c^5, $$
 $$
 be^2f^2, bde^2f, bc^2d^2, b^2e^4, b^2cd^2, b^3de, b^4e^3,
 b^4cf^3, b^5, ae^2f^2, ae^4, ad, $$
 $$
 ac^4e^2, abe^2, abc^4f^2,
 ab^2c^4, ab^3ef^3, a^2f^4, a^2ef^2, a^2c^3, a^3ef, a^3e^3, $$
 $$
 a^3c^2f^2, a^3c^2e, a^3b^3c^2f, a^3b^4c^2, a^4f^3, a^4cf, a^4ce^2,
 a^4b^4c, a^5 \rangle$$
 in $\K[a,b,c,d,e,f]$ defines a $G$-cluster
 $[J]\in \hilbg$.  This can be verified with Macaulay~2~\cite{M2} by
 showing that \(J\) is the initial ideal with respect to the weight
 vector $(22,10,16,50,31,21)$ of the lattice ideal $I_M$.

 To show that \(Y_\theta \cong \hilbg\) is not normal for $\theta\in
 \Theta$ satisfying (\ref{eqn:ghilbtheta}), it is enough by
 Theorem~\ref{thm:coverhilbg} to show that the semigroup \(A_J\) is
 not normal. The set $$\{ (0,-3,0,3,-1,-4), (-3,-3,3,0,-1,4),
 (-4,0,4,1,0,-1),$$
 $$(-2,1,2,2,-3,0), (-1,4,-4,0,3,-2),
 (2,-2,-2,-1,1,2), (3,4,2,-1,-2,-1),$$
 $$(4,-1,-4,0,-2,3), (4,-4,1,-2,2,-1), (3,2,-3,1,-1,-2) \}$$ computed
 using Normaliz generates the semigroup consisting of elements of $M$
 in the rational cone spanned by $A_J$.  Of these vectors, the last
 one $(3,2,-3,1,-1,-2)$ does not lie in $A_J$, hence $A_J$ is not
 normal.
\end{example}

 \begin{remark}
   Example~\ref{ex:nonnormal} was found by applying
   Algorithm~\ref{a:nonnormal}.  The choice of group is a modification
   of an example of a nonnormal toric Hilbert scheme in
   \cite{AlgoTHS}.  We note, however, that the most natural
   modification of that example, using the same weight vector, does
   not work. It is straightforward to modify
   Example~\ref{ex:nonnormal} to get a nonnormal $\hilbg$ for $G
   \subseteq \SL(7,\K)$
 \end{remark}

 This example answers the question of
 Nakamura~\cite[Remark~2.10]{Nakamura}.

 \begin{corollary}
 \label{coro:hilbgnormal}
 Nakamura's \(G\)-Hilbert scheme \(\hilbg\) need not be normal.
 \end{corollary}

 \begin{remark}
 Corollary~\ref{coro:hilbgnormal} implies that the distinguished
 irreducible 
 component $V$ of $Z$ is not normal in general.  This shows that the
 assumption $\N C = \Q_{\geq 0}C\cap \Z C$ made implicitly by Sardo
 Infirri~\cite[Proposition 5.3]{SI2} is not valid in general.
 \end{remark}

 Santos~\cite{Santos} proved that the toric Hilbert scheme may be disconnected. This leads 
naturally to the following conjecture.
                                                                                
 \begin{conjecture}
 There exists a finite abelian subgroup $G\subset \GL(n,\K)$ such that
$\ghilb$ is disconnected.
 \end{conjecture}

\begin{remark}
 For a particular $G \subseteq \GL(n,\K)$ the connectedness of $\ghilb$ can
 be checked by enumerating all monomial ideals on $\ghilb$, and then
 enumerating those in the connected component of the coherent component
 using a modification of the flip graph algorithm from \cite{MT1}.
 Attempting to modify Santos' examples from \cite{Santos} in a similar
 fashion to the above, however, would give a subgroup of $\GL(26,\K)$,
 which is computationally prohibitive to work with.  In addition, just
 as a naive modification of the nonnormal toric Hilbert scheme example does
not give a nonnormal $\hilbg$, there is no reason to expect that this
 subgroup of $\GL(26,\K)$ would have a disconnected $\ghilb$.  The
philosophy remains, however, that multigraded Hilbert schemes tend to be
 disconnected.
\end{remark}

 \end{document}